\documentclass[journal ]{new-aiaa}
\usepackage[utf8]{inputenc}
\usepackage{textcomp}

\usepackage{graphicx}
\usepackage{amsmath}
\usepackage[version=4]{mhchem}
\usepackage{siunitx}
\usepackage{longtable,tabularx}
\usepackage{algpseudocode, algorithm}
\usepackage{multirow, booktabs, makecell}
\usepackage{subfig, pdflscape, float}
\usepackage{booktabs}
\allowdisplaybreaks
\title{Integrated Routing and Trajectory Design for Satellite Servicing}

\author{Euihyeon Choi \footnote{Ph.D. Student, Daniel Guggenheim School of Aerospace Engineering.} and Koki Ho \footnote{Dutton-Ducoffe Professor; Associate Professor, Daniel Guggenheim School of Aerospace Engineering; kokiho@gatech.edu. Associate Fellow AIAA (Corresponding Author).}}
\affil{Georgia Institute of Technology, Atlanta, Georgia 30332}

\begin{document}

\maketitle

\begin{abstract}
This paper studies the integrated spacecraft routing and trajectory optimization problem for satellite servicing missions involving partial en-route propellant replenishment. Unlike terrestrial routing problems, spacecraft operate in a dynamic environment, and we need to optimize the spacecraft routing over a network with nonlinear and time-dependent trajectory costs. In this paper, we tackle this problem using two different formulations. The first formulation, referred to as the arc-based formulation, defines variables based on arcs and employs an iterative decoupling scheme that alternates between mixed-integer linear programming and sequential nonlinear trajectory optimization. The second formulation, referred to as the path-based formulation, defines variables based on paths/routes and leverages column generation and a labeling algorithm to accelerate the identification of promising routes. Through a geosynchronous satellite servicing case study and numerical experiments, we quantify the computational trade-offs between these two formulations in terms of the solution optimality, computational time, and robustness against non-converging or trivial solutions.
\end{abstract}

\section*{Nomenclature}
{\renewcommand\arraystretch{1.0}
\noindent\begin{longtable*}{@{}l @{\quad=\quad} l@{}}
$A_t, \mathbf{a}_{t,\omega}$ & target satellite matrix for feasible routes and the corresponding column $\omega$\\
$A_r, \mathbf{a}_{r,\omega}$ & refueling station matrix for feasible routes and the corresponding column $\omega$\\
$c_\omega$      & corrected profits for route $\omega$\\
$m^{\text{dry}}$ & dry mass of spacecraft \\
$m_{\text{max}}$ & maximum mass of the spacecraft\\
$m_\omega^{*p}$ & initial propellant mass loaded at the starting depot for route $\omega$\\
$n_{dv}, n_{rv}$ & number of duplicated depots and refueling stations\\
$n_r, n_t$     & number of refueling stations and target satellites\\
$p_i$           & obtainable profit from satellite $i$\\
$q_i$           & remaining payload mass when arriving at node $i$\\
$q_{\text{max}}$ & maximum payload capacity (mass) of the spacecraft\\
$R_0, \mathbf{r}_{0,\omega}$ & refueling mass matrix for feasible routes and the corresponding column $\omega$\\
$r_i$           & refueled mass at refueling station $i$\\
$r_{\omega,i}^{*}$ & refueled mass at refueling station $i$ for given route $\omega$\\
$S_{D0}, S_{Dv}, S_{Ds}, S_{De}$ & index sets of original, virtual, starting, and ending depots\\
$S_{R0}, S_{Rv}, S_R$ & index set of original, virtual, and all refueling stations\\
$S_{R,i}$       & index set of duplicated refueling stations corresponding original refueling station $i$\\
$S_T$           & index set of all target satellites \\
$S_{T,\omega}$  & index set of target satellites in route $\omega$ \\
$s_i$           & payload mass for satellite $i$\\
$T^d_i, T^a_i, T^s_i$ & departure, arrival, and servicing time of node $i$\\
$T^t_{ij}$      & orbit transfer time from node $i$ to $j$ \\
$T_{\text{max}}$ & maximum time for the entire mission\\
$u_i$           & spacecraft mass when arriving at node $i$\\
$V$             & index set of all refueling stations and all target satellites.\\
$V'$            & index set of all refueling stations and a subset of target satellites\\
$x_{ij}^k$      & binary decision variable (1 if spacecraft $k$ travels from node $i$ to $j$, and 0 otherwise)\\
$y_{kj}$        & continuous decision variable for linearizing $q_j x_{kj}^k$\\
$z_\omega$      & binary decision variable (1 if route $\omega$ is selected, and 0 otherwise)\\
$\Delta v_{ij}, \mu_{ij}$ & required $\Delta v$ and mass ratio for travel from node $i$ to $j$\\
$\lambda$       & weighted factor of the objective function\\
$\Omega$        & index set of all possible routes\\
\multicolumn{2}{@{}l}{\textit{Subscripts}} \\
$i, j$ & index representing a node \\
$k$ & index representing a depot or corresponding vehicle \\
$m$ & vectors and matrices associated with the restricted master problem \\
$\omega$ & index representing a route (column) \\
\multicolumn{2}{@{}l}{\textit{Superscripts}} \\
$*$ & optimal value \\
\end{longtable*}}
\setcounter{table}{0}

\section{Introduction}
\lettrine{A}{s} the number of satellites expands rapidly for the new space era, the demand for sustainable in-space logistics increases for satellite servicing, repair, and propellant replenishment \cite{Ho_2024_SpaceLogisticsReview, Choi_2025_OrbitalDepotLocation}. Particularly, operational life extension of high-value assets in geostationary orbit or global navigation satellite systems is no longer a laboratory-level concept, but a necessary strategic operation \cite{arney_isam_2024}. However, space missions involve extremely expensive orbital transfers while resources are highly constrained, necessitating meticulously designed high-level mission planning. For instance, utilizing orbital depots and refueling stations for sustainable space missions could be considered \cite{Sarton_2021, Sarton_on-orbit_2022, Choi_2026_VRTPPPR_SciTech}. Such new systems enable missions previously impossible due to spacecraft's mass limitations and can be operated more efficiently in the long term.

For in-space servicing missions, it is essential to identify optimal routes or sequences for spacecraft operations, particularly given resource constraints and limited launch windows \cite{Ahn_2012_GLRPP, Ho_2024_SpaceLogisticsReview}. For example, in practical planetary exploration scenarios that may involve multiple agents and the synergistic benefits of visiting various exploration sites, systematic modeling and solution strategies are necessary \cite{Choi_2022_VRPPCS, Choi_2023_2EVRPP}.

Recent literature has approached the spacecraft routing problem by drawing on methods developed for the terrestrial vehicle routing problem (VRP). The classical VRP aims to minimize overall travel costs while ensuring that all customers are visited, subject to constraints such as maximum capacity or travel distance limits. In the previous literature, various extensions of the VRP have been introduced for practical scenarios and applications \cite{Braekers_2016_VRPReview, Vidal_2020_VRP}. One of the significant variants is the vehicle routing problem with profits (VRPP), which allows abandoning certain tasks (customers) due to limited fuel or resources \cite{Archetti_2014_VRPPCh10}. This characteristic makes the VRPP particularly applicable to space missions, where resource constraints are strict, requiring careful target selection to enhance scientific or economic outcomes. Another important variant addresses the need for refueling during travel, which is extensively explored in operations research for alternative fuel vehicles, known as the green vehicle routing problem (G-VRP) \cite{ERDOGAN_2012_GVRP, ASGHARI_2021_GVRP, MOGHDANI_2021_GVRP}. Initially devised for short-range electric vehicles, the G-VRP's primary feature serves as a useful model for in-space refueling systems, potentially enabling spacecraft to refuel during the mission. 

Despite the close relationship between the spacecraft routing problem and the terrestrial VRP, the unique challenges associated with orbital mechanics make applying conventional VRP approaches to space missions difficult \cite{chen_2021_opt-multitarget-rendz}. Namely, unlike terrestrial environments with fixed distances (costs), in-space operations are dynamic, meaning the cost ($\Delta v$) between nodes is a nonlinear function of departure and transfer times \cite{bang_two-phase_2018, Bang_2019_MultitargetRendezvous}. Consequently, a satellite servicing mission requires the optimization of the visiting sequence over a network with nonlinear and time-dependent trajectories between nodes. Some previous studies have proposed separating these problems by calculating all potential transfers in advance and storing them in a database for VRP \cite{Bang_2019_MultitargetRendezvous, Lee_2023_MultitargetRendezvousHybrid}; however, this approach can be impractical when an intractable number of feasible transfer options exist and cannot be fully enumerated beforehand. Barea et al. \cite{Barea_2020_MultitargetSpaceDebrisRemoval} utilized a two-level approach in which object selection relies on fixed, precomputed $\Delta v$ estimates. While the lower level verifies feasibility, it provides no feedback to update the original cost matrix. So, it might miss efficient trajectories by updating the $\Delta v$ matrix. Recently, Bendimerad et al. \cite{Bendimerad_2025_RefuelingStrategies} proposed an analytic derivation of optimal refueling strategies for an electric propulsion system, focusing on the distribution of a given total $\Delta v$ across $n$ segments for a single transfer to maximize payload and structure mass. However, optimizing both the spacecraft's travel sequence (routes) and trajectories with refuelable infrastructures remains a challenging problem, and we lack an integrated approach that can efficiently optimize them in an integrated way.

To address this gap, this paper formalizes a novel problem, referred to as the vehicle routing and trajectory problem with profits and partial refueling (VRTPP-PR). The VRTPP-PR maximizes total profits and minimizes propellant consumption by determining the selection and order of visited satellites, the refueling amounts at each station, and the departures and transfer times between nodes. This work has three main contributions. First, we propose the VRTPP-PR using two distinct formulations (arc-based and path-based) to capture the integrated nature of spacecraft routing, en-route refueling, and nonlinear trajectory optimization. Second, we develop decomposition-based solution methodologies for each formulation. Finally, through a geosynchronous satellite servicing case study and scaled numerical experiments, we quantify the algorithmic trade-offs between the two approaches.

The novelty of this paper lies not only in proposing a new integrated spacecraft routing and trajectory optimization problem, but also in analyzing the two fundamental types of formulations for this problem: the arc-based and path-based formulations. Many terrestrial VRPs use the path-based formulation \cite{simchi-levi_logic_2014}, but modern space logistics studies have primarily used arc-based formulations \cite{Ho_2024_SpaceLogisticsReview}. While these two formulations have been extensively compared in conventional terrestrial operations research \cite{Masson_2014_1arc2path_PDP, Faiz_2019_cg_vrp}, their comparison in the spacecraft routing problem, which is closely coupled with nonlinear orbital mechanics and in-space infrastructure, remains unexplored. This comparison is nontrivial because it is impacted not only by the routing optimization formulations themselves but also by their treatment of the nonlinear trajectory costs. Because the choice of mathematical formulation significantly impacts computational performance and solution robustness, this paper provides a quantitative comparison of these formulations in the context of the spacecraft routing and trajectory problem as well as a discussion on their pros and cons.

The remainder of this paper is structured as follows. Section \ref{Sec:prob_des} formalizes the VRTPP-PR and presents both arc-based and path-based mathematical formulations. In Sec.~\ref{Sec:sol_meth}, the detailed methodology for each formulation is proposed. Section~\ref{Sec:case_study} presents the results of the case study and numerical experiment, and Sec.~\ref{Sec:conclusions} concludes the whole paper.

\section{Problem Description} \label{Sec:prob_des}

The VRTPP-PR builds on the traditional VRPP by expanding its application to orbital operations, including trajectory optimization and en-route partial refueling for transfers between nodes. The servicer departs from a designated depot and travels to a series of satellites that require service (e.g., refuel or repair). If needed, it can stop at one or more orbiting refueling stations to obtain extra propellant, with the amount of fuel taken as a variable to be optimized. Since greater mass requires more propellant, it is essential to replenish only what is necessary to keep overall fuel consumption to a minimum. Moreover, because each refueling station has limited capacity, mission planners must appropriately determine the amount of refuel at each station. The primary goal of this problem is to maximize total profit from servicing targeted satellites while minimizing fuel consumption, which can be combined into a single objective function using a weighted-sum approach. It is important to note that profits can be generated only through satellite servicing, as this is the key task set by stakeholders. In contrast, refueling stations function solely to replenish fuel and do not generate profit. To clarify, a depot serves as a primary base for initial payload staging and final return, whereas the refueling stations are intermediate nodes for en-route propellant replenishment, analogous to gas stations in terrestrial logistics. In the mathematical formulation, the depot is the start and end node where boundary conditions must be satisfied, while the refueling stations are optional intermediate nodes that extend the spacecraft's travel range.

The cost ($\Delta v$) associated with the trajectory between two nodes (i.e., depot, refueling stations, and satellites) is significantly influenced by the departure and transfer times. To minimize fuel consumption, the spacecraft should strategically wait for a certain time and transfer at the optimal time of flight. Consequently, it is essential to incorporate departure and transfer times into the optimization framework, as the trajectory is closely linked to the sequence of visits.

This study is based on four main assumptions for simplicity: 1) the required payload at each satellite is predetermined, 2) each refueling station has a limit on the number of visits for replenishment, 3) the times for refueling and servicing at each station and satellite are fixed in advance, and 4) transfers are performed using a two-impulse maneuver with a chemical propulsion system. Nonetheless, the proposed frameworks can be easily adapted to relax these assumptions, while the mathematical formulations may need to be adjusted accordingly.

\begin{figure}
    \centering
    \includegraphics[width=3.25in]{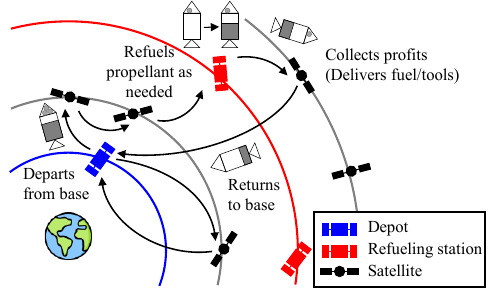}
    \caption{Problem Structure of the VRTPP-PR}
    \label{fig:problem-structure}
\end{figure}

The illustration in Figure \ref{fig:problem-structure} presents a typical example of the VRTPP-PR. This scenario contains a single depot, two refueling stations, and five satellites. The spacecraft begins its mission from the depot, staged with an optimized amount of propellant and the necessary payloads for the multiple services. As depicted in the upper portion of the figure, the spacecraft delivers payloads to the designated satellites, stops at a refueling station to replenish its propellant, then continues to the remaining satellites and returns to the depot. Note that the depot, refueling stations, and satellites in the figure are positioned at non-scale, and the paths between them are shown as arrows for visualization purposes.

\subsection{Arc-based Mathematical Formulation}\label{Sec:ProbDesc_ArcMath}
The arc-based formulation addresses the VRTPP-PR with a static directed graph. To accommodate multiple visits to certain nodes (e.g., the deployment of multiple spacecraft and the re-visitation of refueling stations), we introduce duplicated or virtual depots and refueling stations. Furthermore, we differentiate between the starting and ending depots to accurately calculate the spacecraft's mass and payload history. For clarity, we define the index sets as follows: the original depot set ($S_{D0}$), the virtual depots set ($S_{Dv}$), the starting depots set ($S_{Ds}$), the ending depots set ($S_{De}$), the original refueling stations set ($S_{R0}$), the virtual refueling stations set ($S_{Rv}$), all refueling stations set ($S_R$), the target satellites set ($S_T$), and the collective set of all nodes excluding depots (i.e., all refueling stations and target satellites) ($V$).
\begin{gather}
    S_{D0} = \{0\}\\
    S_{Dv} = \{1, \dots, n_{dv} - 1\}\\
    S_{Ds} = S_{D0} \cup S_{Dv} = \{0, \dots, n_{dv} - 1\}\\
    S_{De} = \{n_{dv}, \dots, 2 n_{dv} - 1\}\\
    S_{R0} = \{2 n_{dv}, \dots, 2 n_{dv}+ n_{r} - 1\}\\
    S_{Rv} = \{2 n_{dv}+ n_{r}, \dots, 2 n_{dv}+ n_{r} n_{rv} - 1\}\\
    S_{R} = S_{R0} \cup S_{Rv} = \{2 n_{dv}, \dots, 2 n_{dv}+ n_{r} n_{rv} - 1\}\\
    S_T = \{2 n_{dv}+ n_{r}n_{rv}, \dots, 2 n_{dv}+ n_{r}n_{rv}+ n_t - 1\}\\
    V = S_{R} \cup S_T = \{2 n_{dv}, \dots, 2 n_{dv}+ n_{r}n_{rv}+ n_t - 1\}
\end{gather}
where $n_{dv}$ and $n_{rv}$ represent the quantity of duplicated depots and refueling stations, including the original depot and refueling stations ($n_{dv} \geq 1$ and $n_{rv} \geq 1$). $n_r$ denotes the count of original refueling stations, and $n_t$ refers to the total number of target satellites. It is important to note that any spacecraft that departs from a starting depot ($k \in S_{Ds}$) is associated with an ending depot ($k'(k) = k + n_{dv}$). As one spacecraft corresponds to a single (duplicated) depot, the number of duplicated depots is the same as the number of available spacecraft, while the duplicated refueling stations indicate the maximum number of potential revisits for each refueling station. For notational convenience, define an index set of all duplicated refueling stations corresponding to the original refueling station $i$ as follows: 
\begin{align}
    S_{R,i} = \{i, i + n_{r}, \dots, i + (n_{rv} - 1)n_{r}\} \quad \forall i \in S_{R0}
\end{align}

The arc-based formulation of the VRTPP-PR can be represented as a mixed-integer nonlinear programming (MINLP), as follows:

($P_a$) Arc-based formulation of the VRTPP-PR
\begin{align}
    \max_{\mathbf{x}, \mathbf{u}, \mathbf{q}, \mathbf{r}, \mathbf{T}^d, \mathbf{T}^t} \sum_{j \in S_T} p_j \left[ \sum_{k \in S_{Ds}} \left(x_{kj}^k + \sum_{i \in V, i \neq j} x_{ij}^k \right)\right] - \lambda \left[ \sum_{k \in S_{Ds}} \left(u_k - m^{\text{dry}} \sum_{j \in V} x_{kj}^k - \sum_{j \in V} q_j x_{kj}^k \right) + \sum_{i \in S_{R}} r_i \right] \label{eq:arc-obj}
\end{align}
subject to
\begin{align}
    & \sum_{j \in V} x_{kj}^k \leq 1 && \forall k \in S_{Ds} \label{eq:arc-const1} \\
    & \sum_{k \in S_{Ds}} \left(x_{kj}^k + \sum_{i \in V, i \neq j} x_{ij}^k \right) \leq 1 && \forall j \in S_{T} \label{eq:arc-const2}\\
    & \sum_{k \in S_{Ds}} \left(x_{kj}^k + \sum_{i \in V, i \neq j} x_{ij}^k \right) \leq 1 && \forall j \in S_{R} \label{eq:arc-const3}\\
    & x_{kj}^k - x_{jk'(k)}^k + \sum_{i \in V, i\neq j} \left(x_{ij}^k - x_{ji}^k\right) = 0 && \forall j \in V, k \in S_{Ds} \label{eq:arc-const4}\\
    & -m_{\text{max}} (1 - x_{kj}^k) \leq u_j - \mu_{kj}(T^d_k, T^t_{kj}) u_k \leq m_{\text{max}} (1 - x_{kj}^k) && \forall j \in V, k \in S_{Ds} \label{eq:arc-const5}\\
    & -m_{\text{max}} \left(1 - \sum_{k \in S_{Ds}} x_{ij}^k \right) \leq u_j - \mu_{ij}(T^d_i, T^t_{ij}) (u_i - s_i) \leq m_{\text{max}} \left(1 - \sum_{k \in S_{Ds}} x_{ij}^k \right) && \forall i \in S_T, j \in V, i \neq j \label{eq:arc-const6}\\
    & -m_{\text{max}} \left(1 - \sum_{k \in S_{Ds}} x_{ij}^k\right) \leq u_j - \mu_{ij}(T^d_i, T^t_{ij}) (u_i + r_i) \leq m_{\text{max}} \left(1 - \sum_{k \in S_{Ds}} x_{ij}^k\right) && \forall i \in S_{R}, j \in V, i \neq j \label{eq:arc-const7}\\
    & -m_{\text{max}} (1 - x_{ik'(k)}^k) \leq m^{\text{dry}} - \mu_{ik'(k)}(T^d_i, T^t_{ik'(k)}) (u_i - s_i) \leq m_{\text{max}} (1 - x_{ik'(k)}^k) && \forall i \in S_T, k \in S_{Ds} \label{eq:arc-const8}\\
    & -m_{\text{max}} (1 - x_{ik'(k)}^k) \leq m^{\text{dry}} - \mu_{ik'(k)}(T^d_i, T^t_{ik'(k)}) (u_i + r_i) \leq m_{\text{max}} (1 - x_{ik'(k)}^k) && \forall i \in S_{R}, k \in S_{Ds} \label{eq:arc-const9}\\
    & u_i \geq m^{\text{dry}} + q_i && \forall i \in {V} \label{eq:arc-const10}\\
    & u_i + r_i \leq m_{\text{max}} && \forall i \in S_{R} \label{eq:arc-const12}\\
    & -q_{\text{max}} \left(1 - \sum_{k \in S_{Ds}} x_{ik'(k)}^k\right) \leq q_i - s_i \leq q_{\text{max}} \left(1 - \sum_{k \in S_{Ds}} x_{ik'(k)}^k\right) && \forall i \in S_T \label{eq:arc-const13}\\
    & -q_{\text{max}} \left(1 - \sum_{k \in S_{Ds}} x_{ik'(k)}^k\right) \leq q_i \leq q_{\text{max}} \left(1 - \sum_{k \in S_{Ds}} x_{ik'(k)}^k\right) && \forall i \in S_{R} \label{eq:arc-const14}\\
    & -q_{\text{max}} \left(1 - \sum_{k \in S_{Ds}} x_{ij}^k\right) \leq q_i - (q_j + s_i) \leq q_{\text{max}} \left(1 - \sum_{k \in S_{Ds}} x_{ij}^k\right) && \forall i \in S_T, j \in V, i \neq j \label{eq:arc-const15}\\
    & -q_{\text{max}} \left(1 - \sum_{k \in S_{Ds}} x_{ij}^k\right) \leq q_i - q_j \leq q_{\text{max}} \left(1 - \sum_{k \in S_{Ds}} x_{ij}^k\right) && \forall i \in S_{R}, j \in V, i \neq j \label{eq:arc-const16}\\
    & \sum_{j \in S_{R,i}} r_{j} \leq r_{\text{max}, i} && \forall i \in S_{R0} \label{eq:arc-const17}\\
    & x_{kj}^k \in \{0, 1\} && \forall j\in V, k \in S_{Ds} \label{eq:arc-const18}\\
    & x_{ij}^k \in \{0, 1\} && \forall i, j\in V, i \neq j, k \in S_{Ds} \label{eq:arc-const19}\\
    & x_{ik'(k)}^k \in \{0, 1\} && \forall i \in V, k\in S_{Ds} \label{eq:arc-const20}\\
    & 0 \leq u_i \leq m_{\text{max}} && \forall i \in S_{Ds} \cup V \label{eq:arc-const21}\\
    & 0 \leq q_i \leq q_{\text{max}} && \forall i \in V \label{eq:arc-const22}\\
    & 0 \leq r_i && \forall i \in S_{R} \label{eq:arc-const23}
\end{align}

The objective function expressed in Eq.~\eqref{eq:arc-obj} aims to maximize the total profits, denoted in the first summation, while concurrently minimizing overall propellant consumption, calculated in the second summation. Decision variables are represented as vectors for each of the variables (e.g., $\mathbf{x} = [x_{ij}^k]$). $p_i$ denotes the profit corresponding to satellite $i$, while $x_{ij}^k$ is a binary variable that indicates whether arc $i$ to $j$ is traveled by spacecraft $k$. The parameter $\lambda$ represents the weighting factor that balances profit maximization against propellant usage; $u_k$ indicates the initial mass at depot $k$; $m^{\text{dry}}$ is the spacecraft's dry mass; $q_j$ specifies the remaining payload mass when arriving node $j$; and $r_i$ quantifies the refueled amount at refueling station $i$.

The constraints in Eqs.~(\ref{eq:arc-const1}--\ref{eq:arc-const3}) ensure that each starting depot, refueling station, and satellite should be visited at most once. Equation \eqref{eq:arc-const4} maintains flow continuity for each node. The constraints in Eqs.~(\ref{eq:arc-const5}--\ref{eq:arc-const7}) enforce the rocket equation for spacecraft transfer between nodes. Here, $u_i$ represents the mass arriving at node $i$, and $\mu_{ij}(T^d_i, T^t_{ij})$ illustrates the mass ratio for transferring from node $i$ to $j$, which varies depending on the departure time from node $i$ ($T^d_i$) and the transfer time to node $j$ ($T^t_{ij}$) \{i.e., $\mu_{ij}(T^d_i, T^t_{ij}) = \exp[-\Delta v_{ij}(T^d_i, T^t_{ij}) / (g_0 I_{sp})]$\}; $m_{\text{max}}$ indicates the spacecraft's maximum mass, and $s_i$ is the payload mass for servicing satellite $i$. 
Note that $\Delta v_{ij}$ is computed by solving Lambert's problem for the specified departure and transfer times. The trajectory design is implicitly embedded in $P_a$ through the optimization of the continuous-time variables ($\mathbf{T}^d$ and $\mathbf{T}^t$). Because explicitly incorporating the full orbital state variables and dynamics into the MINLP would make the formulation excessively complex, $\Delta v_{ij}$ is treated as a nonlinear function of the time variables in $P_a$. The detailed computation and optimization of $\Delta v_{ij}$ are presented in a later section (Sec.~\ref{Sec:CostMtxUpdate}).
Equations \eqref{eq:arc-const8} and \eqref{eq:arc-const9} require that the spacecraft should be empty (without any propellant or payload) upon returning to the depot for optimal transfer efficiency. Equations (\ref{eq:arc-const10}--\ref{eq:arc-const12}) define the upper and lower bounds for spacecraft mass at each refueling station and satellite. Further, Eqs.~(\ref{eq:arc-const13}--\ref{eq:arc-const16}) regulate the remaining payload mass throughout the travel paths, while Eq.~\eqref{eq:arc-const17} sets a cap on the capacity for each refueling station. Lastly, Eqs.~(\ref{eq:arc-const18}--\ref{eq:arc-const23}) specify the conditions for binary and nonnegative continuous decision variables along with their upper limits.

\subsection{Path-based Mathematical Formulation} \label{Sec:ProbDesc_PathMath}
A path-based formulation can also be used to solve the VRTPP-PR. Unless otherwise specified, notations and variables used in this section carry over from the arc-based formulation in Sec.~\ref{Sec:ProbDesc_ArcMath}. 
Let $\Omega$ be the index set of all possible routes for a single spacecraft. Each route $\omega \in \Omega$ represents an ordered sequence of visited nodes: starting from the depot, visiting target satellites and refueling stations, and arriving at an ending depot. Furthermore, let $\Omega_f \subseteq \Omega$ denote the feasible routes in which the ordered sequence of route $\omega \in \Omega$ satisfies all physical constraints (e.g., maximum mass and mission makespan constraint).

Now, we can formulate the VRTPP-PR in a path-based formulation as follows:

($P_p$) Path-based formulation of the VRTPP-PR
\begin{align}
    \max_{\mathbf{z}} \sum_{\omega \in \Omega_f} c_\omega z_\omega = \mathbf{c}^\top \mathbf{z} \label{eq:path-obj}
\end{align}
subject to
\begin{align}
    &\sum_{\omega \in \Omega_f} \mathbf{a}_{t,\omega} z_\omega = A_{t} \mathbf{z} \leq \mathbf{1}_{n_t} \label{eq:path-const1}\\
    &\sum_{\omega \in \Omega_f} \mathbf{a}_{r,\omega} z_\omega = A_{r} \mathbf{z} \leq n_{rv} \mathbf{1}_{n_r} \label{eq:path-const2}\\
    &\sum_{\omega \in \Omega_f} z_\omega = \mathbf{1}^\top \mathbf{z}  \leq n_{dv} \label{eq:path-const3}\\
    &\sum_{\omega \in \Omega_f} \mathbf{r}_{0,\omega} z_\omega = R_0 \mathbf{z} \leq \mathbf{r}_{\text{max}} \label{eq:path-const4}\\
    &\mathbf{z} \in \{0, 1\}^{|\Omega_f|} \label{eq:path-const5}
\end{align}

Equation~\eqref{eq:path-obj} indicates the maximizing corrected profits (weighted sum of profits and propellant consumption) of all used routes. For a given feasible route $\omega \in \Omega_f$, we define $S_{T, \omega} \subseteq S_T$ as the subset of target satellites visited in the route $\omega$. The definition of the corrected profit of route $\omega$ is as follows:
\begin{align}
    c_{\omega} = \sum_{i \in S_{T,\omega}} p_i - \lambda \left(m^{*p}_\omega + \sum_{i \in S_{R}} r_{\omega,i}^* \right) \label{eq:cw}
\end{align}
where $m^{*p}_\omega$ is the initial propellant mass loaded at the depot and $r^*_{\omega,i}$ is the refueling amount at refueling station $i \in S_R$ required to complete the route $\omega$. $\mathbf{z} = [z_\omega]$ is the binary decision variable vector, while $z_\omega = 1$ if route $\omega$ is used and 0 otherwise. 

Equation~\eqref{eq:path-const1} ensures that each satellite is visited at most once. Here, $\mathbf{a}_{t,\omega} \in \{0, 1\}^{n_t}$ is the incidence vector of the set $S_{T, \omega}$ with respect to $S_T$. That is, the component corresponding to satellite $i \in S_T$ is 1 if $i \in S_{T, \omega}$ and 0 otherwise. Similarly, for Eq.~\eqref{eq:path-const2}, $\mathbf{a}_{r,\omega} \in \mathbb{Z}_{\geq 0}^{n_r}$ represents the visit counts for original refueling stations. The component corresponding to the original refueling station $i$ ($a_{r,\omega,i}$) equals the number of visits to $i$ (including its virtual copies) in route $\omega$. Equation~\eqref{eq:path-const3} restricts the maximum number of spacecraft (routes) of the whole mission, and Eq.~\eqref{eq:path-const4} indicates the maximum capacity of each refueling station. $\mathbf{r}_{0,\omega} \in \mathbb{R}^{n_r}_{\geq 0}$ shows the refueling mass per original refueling station, where its component corresponding to station $i \in S_{R0}$ is the sum of refueling mass of all duplicated refueling stations ($r_{0,\omega, i} = \sum_{j \in S_{R, i}} r^*_{\omega,j}$). $\mathbf{r}_{\text{max}}$ indicates the maximum capacity for each original refuel station. Finally, Eq.~\eqref{eq:path-const5} indicates that the decision variables are binary. 

\section{Solution Methodology} \label{Sec:sol_meth}

The VRTPP-PR simultaneously optimizes visit sequences, transfer trajectories, and refueling strategies in satellite servicing missions. As explained in Sections \ref{Sec:ProbDesc_ArcMath} and \ref{Sec:ProbDesc_PathMath}, the problem can be formulated through either arc-based or path-based methods. However, as the problem itself is a MINLP, the coupled relationships between discrete decisions (such as the visiting sequence) and continuous variables (such as trajectory and refueling amount) make it computationally impractical to solve these models simultaneously due to significant nonlinearity.

To tackle this issue, we propose methodologies devised for each formulation. For the arc-based method, we divide the problem into two distinct subproblems: a mixed-integer linear programming (MILP) that optimizes sequencing and refueling, and a nonlinear programming (NLP) that finds the optimal trajectories between nodes. The two subproblems are solved iteratively until a converged solution is obtained. On the other hand, for the path-based formulation, we adopt a column generation method to manage the exponential number of potential routes. Within this framework, a labeling algorithm is developed to address the route-generation pricing problem, while an NLP solver is employed to refine trajectories, similar to those used in the arc-based method.

\subsection{Solution Methodology for the Arc-based Formulation}
\subsubsection{Initial Cost Matrix Estimation}
To initiate the iterative process, an initial estimate of the cost (mass ratio or $\Delta v$) matrix is essential due to a critical interdependence: the arc-based formulation (MILP) relies on mass ratios to find the optimal sequence of visits, while the NLP necessitates this sequence to refine the trajectory optimizations. To derive this initial matrix, we must evaluate the transfers between each node pair, including depots, refueling stations, and satellites. Given that the actual orbits are generally neither coplanar nor circular, the ideal Hohmann transfer $\Delta v$ cannot be used directly. Consequently, we solve the trajectory optimization problem for each node pair to identify the optimal departure and transfer times, using initial guesses of zero departure time and the Hohmann transfer time.

\subsubsection{Mixed-Integer Linear Programming with Constant Cost Matrix} \label{Sec:MILP}
When we fix the mass ratio matrix, the $P_a$ defined in Section \ref{Sec:ProbDesc_ArcMath} becomes the MILP as follows: 

($P_l$) Linearized formulation of $P_a$
\begin{align}
    \max_{\mathbf{x}, \mathbf{u}, \mathbf{q}, \mathbf{r}, \mathbf{y}} \sum_{j \in S_T} p_j \left[ \sum_{k \in S_{Ds}} \left(x_{kj}^k + \sum_{i \in V, i \neq j} x_{ij}^k \right)\right] - \lambda \left[ \sum_{k \in S_{Ds}} \left(u_k - m^{\text{dry}} \sum_{j \in V} x_{kj}^k - \sum_{j \in V} y_{kj} \right) + \sum_{i \in S_{R}} r_i \right] \label{eq:arc-lin-obj}
\end{align}
subject to
\begin{align}
    &\text{Equations (\ref{eq:arc-const1}--\ref{eq:arc-const4}), (\ref{eq:arc-const10}--\ref{eq:arc-const23})} \nonumber\\
    & -m_{\text{max}} (1 - x_{kj}^k) \leq u_j - \hat{\mu}_{kj} u_k \leq m_{\text{max}} (1 - x_{kj}^k) && \forall j \in V, k \in S_{Ds} \label{eq:arc-lin-const5}\\
    & -m_{\text{max}} \left(1 - \sum_{k \in S_{Ds}} x_{ij}^k \right) \leq u_j - \hat{\mu}_{ij} (u_i - s_i) \leq m_{\text{max}} \left(1 - \sum_{k \in S_{Ds}} x_{ij}^k \right) && \forall i \in S_T, j \in V, i \neq j \label{eq:arc-lin-const6}\\
    & -m_{\text{max}} \left(1 - \sum_{k \in S_{Ds}} x_{ij}^k\right) \leq u_j - \hat{\mu}_{ij} (u_i + r_i) \leq m_{\text{max}} \left(1 - \sum_{k \in S_{Ds}} x_{ij}^k\right) && \forall i \in S_{R}, j \in V, i \neq j \label{eq:arc-lin-const7}\\
    & -m_{\text{max}} (1 - x_{ik'(k)}^k) \leq m^{\text{dry}} - \hat{\mu}_{ik'(k)} (u_i - s_i) \leq m_{\text{max}} (1 - x_{ik'(k)}^k) && \forall i \in S_T, k \in S_{Ds} \label{eq:arc-lin-const8}\\
    & -m_{\text{max}} (1 - x_{ik'(k)}^k) \leq m^{\text{dry}} - \hat{\mu}_{ik'(k)} (u_i + r_i) \leq m_{\text{max}} (1 - x_{ik'(k)}^k) && \forall i \in S_{R}, k \in S_{Ds} \label{eq:arc-lin-const9}\\
    &0 \leq y_{kj} \leq q_{j} && \forall j \in V, k \in S_{Ds} \label{eq:y1}\\
    &q_j - q_{\text{max}}(1 - x_{kj}^k) \leq y_{kj} \leq q_{\text{max}} x_{kj}^k && \forall j \in V, k \in S_{Ds} \label{eq:y2}
\end{align}

In Eq.~\eqref{eq:arc-lin-obj}, additional continuous decision variables $y_{kj}$ are introduced to linearize the multiplications of continuous and binary decision variables in the original objective function [$y_{kj} = q_j x_{kj}^k$ in Eq.~\eqref{eq:arc-obj}]. By using auxiliary linear constraints [Eqs.~(\ref{eq:y1}--\ref{eq:y2})], $y_{kj} = q_j x_{kj}^k$ relationship holds. $\hat{\mu}_{kj}$, $\hat{\mu}_{ij}$, and $\hat{\mu}_{ik'(k)}$ in Eqs.~(\ref{eq:arc-lin-const5}--\ref{eq:arc-lin-const9}), are the constant transfer mass ratios as the departure and transfer times are fixed.

\subsubsection{Cost Matrix Update with Routing Solution} \label{Sec:CostMtxUpdate}
After obtaining a visiting sequence from $P_l$, the next step involves resolving the trajectory optimization problems corresponding to that sequence. To optimize mass ratios, a NLP is solved to minimize $\Delta v$ (or maximize the mass ratio) for each route leg by determining the optimal departure and transfer times. While it is possible to design the NLP to optimize all trajectories along the route concurrently, we take a sequential approach that tackles each segment individually, resulting in faster, more stable convergence. For instance, consider a basic visiting sequence: departing from the depot, visiting a satellite, and returning to the depot. Initially, the departure time from the depot and the transfer time to the satellite are selected to minimize $\Delta v$. Subsequently, these values are fixed, and the departure time from the satellite and the transfer time back to the depot are determined. Note that there is a lower bound on the departure time for subsequent optimization: it must occur after the arrival time and the refueling/servicing time at that node.

The NLP that optimizes the trajectory from node $i$ to $j$ is formulated as follows:

($P_t$) Trajectory optimization problem between two nodes
\begin{align}
    \min_{T^d_i, T^t_{ij}} \Delta v_{ij}(T^d_i, T^t_{ij}) \label{eq:NLP-obj}
\end{align}
subject to 
\begin{gather}
    T^a_i + T^s_i \leq T^d_i \leq T_{\text{max}}\label{eq:NLP-const1}\\
    \epsilon \leq T^t_{ij} \leq T_{\text{max}}\label{eq:NLP-const2}\\
    \epsilon \leq T^d_i + T^t_{ij} \leq T_{\text{max}} - T^s_j\label{eq:NLP-const3}
\end{gather}

$\Delta v$ in the objective function [Eq.~\eqref{eq:NLP-obj}] can be calculated by solving Lambert's problem with specified departure and transfer times \cite{dario_izzo_2020_4091753}. To solve Lambert's problem, it is necessary to determine the position vectors of node $i$ at time $T^d_i$ and of node $j$ at the time $T^d_i + T^t_{ij}$. Given that the orbital elements of all nodes are known in advance, these position vectors can be derived by solving Kepler's problem for the corresponding times \cite{dario_izzo_2020_4091753}. Moreover, to achieve a minimal $\Delta v$, it is assumed that the transfer orbit is in a prograde direction. In Eq.~\eqref{eq:NLP-const1}, $T^a_i$ represents the arrival time at node $i$, which is determined by the solution from the prior trajectory leg, while $T^s_i$ denotes the time for servicing or refueling at node $i$, and $T_{\text{max}}$ refers to the total mission duration. Once the optimal departure time ($T^{d*}_i$) and transfer times ($T^{t*}_{ij}$) for the current segment are identified, the arrival time at the next node $j$ is adjusted for the following optimization step:
\begin{align}
    T^a_j = T^{d*}_i + T^{t*}_{ij}
\end{align}

Equation \eqref{eq:NLP-const2} constrains a lower limit on the transfer time, with $\epsilon$ representing a small positive constant (e.g., $\epsilon = 10^{-5}$ in this paper) to avoid singularities in Lambert's problem. In addition, Eq.~\eqref{eq:NLP-const3} enforces that the sum of the arrival time and servicing/refueling duration at node $j$ must not exceed the mission makespan. The initial estimates for departure and transfer times in the NLP are based on values from the prior iteration. It is noteworthy that, depending on the solver, the solution may not be the global minimum $\Delta v$ trajectory due to the lower/upper bound. In other words, the solver may find the trajectory that can travel faster but consume additional propellant, which is still a local minimum of the constrained NLP. In this study, we utilized a sequential least squares programming (SLSQP) method to solve the NLP, although other suitable nonlinear optimization algorithms could also be used.

\subsubsection{Overall Solution Framework of Arc-based Formulation} \label{Sec:sol-frame-arc}
The mass ratio matrix in the MILP ($P_l$) is updated by solving the sequential $P_t$ for the next iteration. The iterative process continues until convergence. The algorithm is considered to have converged if the visiting sequences are the same and the total change in the $\Delta v$ matrix is below a predefined tolerance threshold, $\epsilon_c$ (e.g., 0.01 km/s), or if a solution cycle is detected. 

\begin{figure}[!hbt]
    \centering
    \includegraphics[width=3.25in]{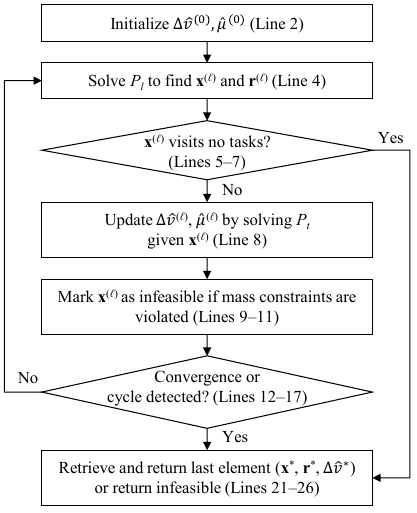}
    \caption{Flow chart of the Algorithm \ref{Alg:arc-based}}
    \label{fig:flow-chart-alg1}
\end{figure}

Algorithm \ref{Alg:arc-based} shows the pseudo-code of the proposed solution methodology for the arc-based formulation and Fig.~\ref{fig:flow-chart-alg1} visualizes the overall procedure. From initialization of $\Delta \hat{v}$ and $\hat{\mu}$ matrices, the linearized arc-based formulation ($P_l$) is solved. Then, for a given visiting sequence, the trajectory optimization problem ($P_t$) is sequentially solved to update $\Delta \hat{v}$ and $\hat{\mu}$ matrices. If there is no visited task, the framework is terminated. Otherwise, the final solution is obtained either if it is the same as the previous solution within a certain threshold, or if the cycle is detected (i.e., if it is the same as one of the past solutions within a certain threshold).
Note that Algorithm 1 does not explicitly generate no-good cuts to prevent the reselection of an infeasible route. Instead, it relies on the iterative $\Delta v$ matrix (equivalently, $\hat{\mu}$ matrix) update. If a refined trajectory causes the spacecraft to exceed its maximum mass limit, the updated (tighter) mass ratio matrix is fed back into the MILP in the next iteration. Because the MILP inherently satisfies the maximum mass constraint, the infeasible route is not selected in the subsequent iteration.

Furthermore, lines \ref{line:cycle1}--\ref{line:cycle2} in Algorithm \ref{Alg:arc-based} implement a cycle-detection mechanism. Because the solution framework decouples the problem into discrete routing and trajectory optimization problems, the fixed-point iteration can occasionally oscillate between two or more suboptimal solutions. To prevent infinite loops, the algorithm tracks the solution history. If the solver revisits a previously evaluated solution (with the same routes and cost matrix), the iteration terminates and returns the most recent solution. While this cycle-breaking logic is heuristic, it is necessary to find a good feasible solution within practical time limits, since the original problem is highly nonlinear and non-convex, making global optimality difficult to ensure.

\begin{algorithm}[!hbt]
\caption{Solution Framework for the Arc-based Formulation}
\label{Alg:arc-based}
\renewcommand{\algorithmicrequire}{\textbf{Input:}}
\renewcommand{\algorithmicensure}{\textbf{Output:}}
\begin{algorithmic}[1]
    \Require Orbital elements of all nodes, problem parameters
    \Ensure Optimal visiting sequence $\mathbf{x}^*$, refueling plan $\mathbf{r}^*$, and trajectory $\Delta \hat{v}^*$
    \State Initialize solution history $\mathcal{H} \leftarrow \emptyset$, iteration counter $\ell \leftarrow 0$
    \State Initialize $\Delta \hat{v}^{(0)}$ and $\hat{\mu}^{(0)}$
    
    \While{$\ell < L_{\text{max}}$}
        \State Solve $P_l$ to determine visiting sequence $\mathbf{x}^{(\ell)}$ and refueling plan $\mathbf{r}^{(\ell)}$ \label{line:Pl} \Comment{MILP routing}
        
        \If{$\mathbf{x}^{(\ell)}$ visits no tasks}
            \State \textbf{break}
        \EndIf
        
        \State Update $\Delta \hat{v}^{(\ell)}$ and $\hat{\mu}^{(\ell)}$ by solving $P_t$ given $\mathbf{x}^{(\ell)}$ \Comment{NLP trajectory refinement}
        
        \If{updated trajectory violates mass constraints}
             \State Mark $\mathbf{x}^{(\ell)}$ as infeasible \Comment{Flag to prevent returning invalid solution}
        \EndIf
        
        \If{$\ell > 0$} 
            \If{$\exists h < \ell$ s.t. $\mathbf{x}^{(\ell)} = \mathbf{x}^{(h)}$ \textbf{and} $\|\Delta \hat{v}^{(\ell)} - \Delta \hat{v}^{(h)}\| < \epsilon_c$} \label{line:cycle1}
                \State Update $\mathcal{H} \leftarrow \mathcal{H} \cup \{(\mathbf{x}^{(\ell)}, \mathbf{r}^{(\ell)}, \Delta \hat{v}^{(\ell)})\}$
                \State \textbf{break} \Comment{Convergence ($h = \ell-1$) or cycle detected ($h < \ell-1$)}
            \EndIf \label{line:cycle2}
        \EndIf 
        
        \State Update $\mathcal{H} \leftarrow \mathcal{H} \cup \{(\mathbf{x}^{(\ell)}, \mathbf{r}^{(\ell)}, \Delta \hat{v}^{(\ell)})\}$
        \State $\ell \leftarrow \ell + 1$
    \EndWhile

    \State Retrieve last element $(\mathbf{x}^*, \mathbf{r}^*, \Delta \hat{v}^*)$ from $\mathcal{H}$
    \If{$\mathbf{x}^*$ is feasible}
        \State \textbf{return} $\mathbf{x}^*$, $\mathbf{r}^*$, $\Delta \hat{v}^*$
    \Else
        \State \textbf{return} Infeasible
    \EndIf
\end{algorithmic}
\end{algorithm}

\subsection{Solution Methodology for the Path-based Formulation}

\subsubsection{Restricted Master Problem Formulation for Column Generation}
The path-based formulation ($P_p$) inherently contains many variables (columns), whereas only a small number are needed to obtain an optimal solution. Therefore, a column generation approach can be used to reduce the problem size, thereby improving computational time. However, since the column generation method is for linear programming (LP), we need to formulate the LP relaxation of $P_p$ as follows:

($P_e$) LP relaxation of $P_p$
\begin{align}
    \max_{\mathbf{z}} \mathbf{c}^\top \mathbf{z} \label{eq:relaxed-obj}
\end{align}
subject to
\begin{align}
    &\text{Equations~(\ref{eq:path-const1}--\ref{eq:path-const4})} \nonumber\\
    &\mathbf{z} \geq \mathbf{0} \label{eq:relaxed-const-1}
\end{align}

Note that the only difference between $P_e$ and $P_p$ is that the binary constraint [Eq.~\eqref{eq:path-const5}] is changed into the non-negativity constraint [Eq.~\eqref{eq:relaxed-const-1}].

Define a restricted master problem (fractional $P_e$) as follows:

($P_m$) Restricted master problem of $P_e$
\begin{align}
    \max_{\mathbf{z}_m} \mathbf{c}^\top_m \mathbf{z}_m \label{eq:master-obj}
\end{align}
subject to
\begin{align}
    A_{t,m} \mathbf{z}_m &\leq \mathbf{1}_{n_t} \label{eq:master-const1}\\
    A_{r,m} \mathbf{z}_m &\leq n_{rv} \mathbf{1}_{n_r} \label{eq:master-const2}\\
    \mathbf{1}^\top_m \mathbf{z}_m  &\leq n_{dv} \label{eq:master-const3}\\
    R_{0,m} \mathbf{z}_m &\leq \mathbf{r}_{\text{max}} \label{eq:master-const4}\\
    \mathbf{z}_m &\geq 0 \label{eq:master-const5}
\end{align}
where the subscript $m$ represents the vectors/matrices that contain a subset of columns of the original formulation in $P_e$. Starting from the initialized columns, the column generation method iteratively selects the most promising column and adds it to the restricted master problem by solving the subproblem described in the next section.

\subsubsection{Initial Column Generation}
The initial columns are generated by solving the vehicle routing and trajectory problem to visit a single satellite $i_0 \in S_T$ using a single spacecraft. This problem can be formulated using the same variables and constraints defined in Section~\ref{Sec:ProbDesc_ArcMath}, with the four adjustments: 1) the objective function is replaced to minimize the overall fuel consumption, 2) the number of duplicated depots is unity ($n_{dv} = 1$) to model a single spacecraft, 3) $S_T$ and $V$ are replaced by the subsets $\{i_0\}$ and $V' = S_R \cup \{i_0\}$, and 4) Eq.~\eqref{eq:arc-const2} is tightened from an inequality into an equality constraint. Mathematically, the subproblem is formulated as follows:

($P_r$) Restricted problem for visiting a single satellite with a spacecraft
\begin{align}
    \min_{\mathbf{x}, \mathbf{u}, \mathbf{q}, \mathbf{r}, \mathbf{T}^d, \mathbf{T}^t} \quad & \sum_{k \in S_{Ds}} \left(u_k - m^{\text{dry}} \sum_{j \in V'} x_{kj}^k - \sum_{j \in V'} q_j x_{kj}^k \right) + \sum_{i \in S_{R}} r_i \label{eq:singe-path-obj}
\end{align}
subject to
\begin{align}
    & \sum_{k \in S_{Ds}} \left(x_{k i_0}^k + \sum_{i \in V', i \neq i_0} x_{i i_0}^k \right) = 1 \label{eq:single-path-const1}\\
    &\text{Equations~\eqref{eq:arc-const1} and (\ref{eq:arc-const3}--\ref{eq:arc-const23}) with $S_T \rightarrow \{i_0\}$ and $V \rightarrow V'$} \nonumber
\end{align}

Equation~\eqref{eq:singe-path-obj} minimizes the total propellant consumption, instead of the profit-maximizing objective, as the designated satellite must be visited in $P_r$. Equation~\eqref{eq:single-path-const1} replaces Eq.~\eqref{eq:arc-const2} with an equality constraint, ensuring that satellite $i_0$ is visited exactly once. All other constraints, such as flow continuity, mass dynamics, payload capacity, and variable bounds [Eqs.~\eqref{eq:arc-const1} and (\ref{eq:arc-const3}--\ref{eq:arc-const23})] remain the same, while the index sets $S_T$ and $V$ are replaced by $\{i_0\}$ and $V'$, respectively.

As with the arc-based method (Sec.~\ref{Sec:sol-frame-arc}), the trajectories must be refined for each generated initial column. Algorithm \ref{Alg:arc-based} is repeated, replacing $P_l$ in Line~\ref{line:Pl} with $P_r$, while the rest of the logic remains the same.

\subsubsection{Subproblem Formulation for Column Generation}
The column generation uses the primal-dual relationship. If there exists a column violating the dual problem's constraint, that column must be included in the primal problem to further increase the objective value in the maximization problem \cite{bertsimas-LPbook, simchi-levi_logic_2014}. Consider the dual problem of $P_e$ as follows:

($P_d$) The dual problem of $P_e$ 
\begin{align}
    \min_{\mathbf{w}_1, \mathbf{w}_2, w_3, \mathbf{w}_4} \mathbf{w}_1^\top \mathbf{1}_{n_t} + n_{rv} \mathbf{w}_2^\top \mathbf{1}_{n_{r}} + w_3 n_{dv} + \mathbf{w}_4^\top \mathbf{r}_{\text{max}}
\end{align}
subject to
\begin{align}
    &A_{t}^\top \mathbf{w}_1 + A_{r}^\top \mathbf{w}_2 + w_3 \mathbf{1} + R_0^\top \mathbf{w}_4 \geq \mathbf{c}\\
    &\mathbf{w}_1, \mathbf{w}_2, w_3, \mathbf{w}_4 \geq 0
\end{align}
where $\mathbf{w}_1$, $\mathbf{w}_2$, $w_3$, and $\mathbf{w}_4$ are the dual variables corresponding to the primal problem constraints [Eqs.~(\ref{eq:path-const1}--\ref{eq:path-const4})]. The current optimal dual variables ($\mathbf{w}_1^*$, $\mathbf{w}_2^*$, $w_3^*$, and $\mathbf{w}_4^*$) can be obtained by solving $P_m$ and a new column (route) $\omega$ satisfying the following condition should be generated:
\begin{align}
    \mathbf{a}_{t,\omega}^\top \mathbf{w}_1^{*} + \mathbf{a}_{r,\omega}^\top \mathbf{w}_2^{*} + w_3^* + \mathbf{r}_{0,\omega}^\top \mathbf{w}_4^{*}  < c_\omega \label{eq:cg-condition}
\end{align}

Among the columns satisfying Eq.~\eqref{eq:cg-condition}, the most violating column (highest reduced cost) $\omega$ is added for each iteration. Note that since $P_e$ is a maximization problem, this violating quantity technically represents a marginal profit, which aims to maximize. However, in this paper, we use the term \textit{reduced cost} to retain standard linear programming terminology. Since one can express $\mathbf{a}_{t,\omega}^\top \mathbf{w}_1^{*}$, $\mathbf{a}_{r,\omega}^\top \mathbf{w}_2^{*}$, and $\mathbf{r}_{0,\omega}^\top \mathbf{w}_4^{*}$ as:
\begin{align}
    \mathbf{a}_{t,\omega}^\top \mathbf{w}_1^{*} &= \sum_{i \in S_{T,\omega}} w_{1,i}^*\\
    \mathbf{a}_{r,\omega}^\top \mathbf{w}_2^{*} &= \sum_{i \in S_{R0}} a_{r,\omega,i} w_{2, i}^* \label{eq:w2-ar}\\
    \mathbf{r}_{0,\omega}^\top \mathbf{w}_4^{*} &= \sum_{i \in S_{R0}} r_{0,\omega,i} w_{4,i}^* \label{eq:w4-rw}
\end{align}
where $w_{1,i}^*$, $w_{2,i}^*$, and $w_{4,i}^*$ are the $i$th component of dual variable vectors $\mathbf{w}_1^{*}$, $\mathbf{w}_2^{*}$, and $\mathbf{w}_4^{*}$, respectively. With Eq.~\eqref{eq:cw} in consideration, the column generation problem (subproblem) can be formulated as follows:

($P_s$) Subproblem to generate column
\begin{align}
    \max_{\omega \in \Omega_f} \sum_{i \in S_{T,\omega}} (p_i - w_{1,i}^*) - \sum_{i \in S_{R0}} (a_{r,\omega,i} w_{2, i}^* + r_{0,\omega,i} (\lambda + w_{4,i}^*)) - w_3^* - \lambda m_\omega^{*p} \label{eq:reduced-cost}
\end{align}

$P_s$ is an elementary shortest-path problem with resource constraints. In this paper, a backward labeling algorithm is used to find the maximum reduced cost column (Algorithm \ref{Alg:labeling}). The backward approach constructs the path from the ending depot (index: 1) back to the starting depot (index: 0). The backward direction is computationally advantageous because the spacecraft mass at the destination/refueling station is fixed at the dry mass ($m^{\text{dry}}$), whereas the departure mass is a variable dependent on the total path fuel consumption. 

Algorithm \ref{Alg:labeling} is designed to find the column with the maximum reduced cost. A straightforward approach to maximize Eq.~\eqref{eq:reduced-cost} is to compute and track the reduced cost directly at each intermediate node along the path. However, we cannot estimate the total initial propellant $m_\omega^{*p}$ exactly during a backward search, as it depends on the remaining (unvisited) sequence of nodes. Therefore, we decompose the objective and define an accumulated reward ($\phi$) from node $i$ to the ending depot, which excludes the $-\lambda m_\omega^{*p}$ term in Eq.~\eqref{eq:reduced-cost}:
\begin{align}
    \phi = \sum_{j \in \text{set}(\boldsymbol{\gamma}) \cap S_T} (p_j - w_{1,j}^*) - \sum_{j \in \text{set}(\boldsymbol{\gamma}) \cap S_{R}} \left(w_{2,{j'(j)}}^* + r_{j}(\lambda +  w_{4,{j'(j)}}^*)\right) - w_3^* \label{eq:phi_i}
\end{align}
where $\boldsymbol{\gamma}$ is the visited node sequence from the ending depot to node $i$ and $j'(j) = \bmod(j - 2 n_{dv}, n_r) + 2 n_{dv}$ indicates the original refueling station index corresponding to the refueling station index $j \in S_R$. The goal of the Algorithm \ref{Alg:labeling} is to maximize the final reduced cost, expressed as $\bar{c} = \phi - \lambda m_\omega^{*p}$ at the starting depot. To eventually calculate $m_\omega^{*p}$, we track the arriving spacecraft mass ($u$) at node $i$ along the visiting sequence. When the backtracking algorithm finally arrives at the starting depot, $m_\omega^{*p}$ can be determined by subtracting the dry mass and cumulative payload mass from $u$.

\begin{figure}[!hbt]
    \centering
    \includegraphics[width=3.25in]{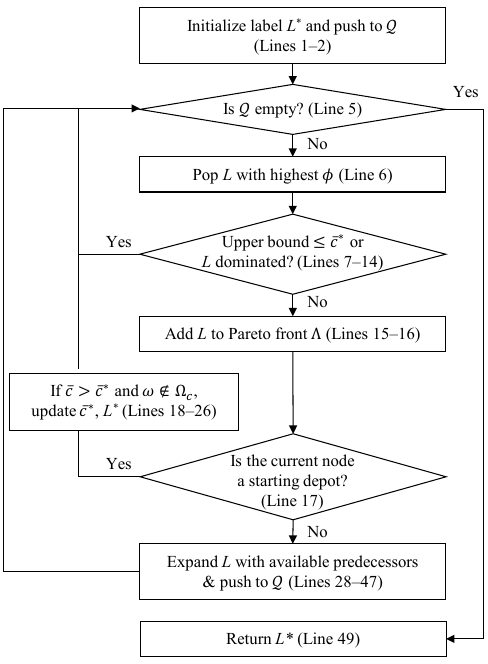}
    \caption{Flow chart of the Algorithm \ref{Alg:labeling}}
    \label{fig:flow-chart-alg2}
\end{figure}

Similar to a tree-search algorithm, Algorithm \ref{Alg:labeling} starts from the ending depot and expands all node candidates until it reaches the starting depot.  Figure~\ref{fig:flow-chart-alg2} illustrates the overall backward search process. To track necessary quantities, we define a label $L$ as a tuple of current node $i$, accumulated reward ($\phi$), arriving spacecraft mass ($u$), path vector ($\boldsymbol{\gamma}$), and the refueling amount vector ($\boldsymbol{\rho}$): 
\begin{align}
    L = (i, \phi, u, \boldsymbol{\gamma}, \boldsymbol{\rho})
\end{align}

Initially, the label starts from the ending depot with the spacecraft's dry mass ($u = m^{\text{dry}}$) and initial reward offset $(\phi = - w_3^*)$, as indicated in Line~\ref{line:initial} and moves backward. Let the algorithm be at node $i$ with label $L$ and be investigating the predecessor node $j$. If the node $j$ is a target satellite, $(p_j - w_{1,j}^*)$ is added to $\phi$ (Line~\ref{line:dual profit}), and $u$ is updated to $u/\hat{\mu}_{ji} + s_j$ (Lines~\ref{line:u-update} and \ref{line:uj target}). Conversely, if node $j$ is a refueling station, $\phi$ is decreased by the refueling-related costs, $(w_{2,{j'(j)}}^* + r_{j}(\lambda + w_{4,{j'(j)}}^*))$ (Line~\ref{line:refueling-related cost}), and $u$ becomes $u/\hat{\mu}_{ji} - r_j$ (Lines~\ref{line:u-update} and \ref{line:uj refuel}). Here, $r_j$ is the refueling amount at a station $j \in S_R$, determined by the available station capacity and the spacecraft's fuel deficit:
\begin{align}
    r_j = \min \left\{\left(r_{\text{max}, j'(j)} - \sum_{\ell \in S_{R,j'(j)}} \rho_{\ell}\right), \left(\frac{u}{\hat{\mu}_{ji}}- m^{\text{dry}} - \sum_{{\ell} \in \text{set}(\boldsymbol{\gamma}) \cap S_T} s_\ell\right)\right\} \label{eq:refueling-rule}
\end{align}
where the first term represents the remaining available refueling amount, and the second term indicates the spacecraft must deplete all the propellant when arriving at the refueling station for efficiency.

Let the label $L_1$ at node $i$ with path $\boldsymbol{\gamma}_1$ be currently under consideration. There might exist another label $L_2$ at node $i$ with a different path sequence ($\boldsymbol{\gamma}_2 \neq \boldsymbol{\gamma}_1$) but the same visited nodes [$\text{set}(\boldsymbol{\gamma}_1) = \text{set}(\boldsymbol{\gamma}_2)$]. Note that any difference in $\phi$ is strictly due to the refueling cost [the second term in Eq.~\eqref{eq:phi_i}] because the collected target rewards [the first term in Eq.~\eqref{eq:phi_i}] are identical for both labels. Consider the case where $u_1 \geq u_2$. If the spacecraft travels the same remaining path from node $i$ to the starting depot, departing with a heavier mass $u_1$ will require a larger final $m_\omega^{*p}$ or greater intermediate refueling than departing with $u_2$ since $r_j$ is deterministic from Eq.~\eqref{eq:refueling-rule}. Since both $m_\omega^{*p}$ and refueling amount are penalty terms in calculating the reduced cost [Eq.~\eqref{eq:reduced-cost}], if $\phi_1 \leq \phi_2$ and $u_1 \geq u_2$, label 2 dominates label 1, so we can prune label 1 (Line \ref{line:dominate check}). On the other hand, if ($\phi_1 > \phi_2$ and $u_1 \geq u_2$) or ($\phi_1 \leq \phi_2$ and $u_1 < u_2$), there is a trade-off; both labels are non-dominated. This Pareto front maintenance accelerates the algorithm by pruning.

\begin{algorithm}[htb!]
\caption{Backward Labeling Algorithm For Column Generation}\label{Alg:labeling}
\renewcommand{\algorithmicrequire}{\textbf{Input:}}
\renewcommand{\algorithmicensure}{\textbf{Output:}}
\small
\begin{algorithmic}[1]
\Require Dual variables ($\mathbf{w}_{1}^*, \mathbf{w}_{2}^*, w_{3}^*, \mathbf{w}_{4}^*$) and previously calculated routes ($\Omega_c$)
\Ensure Optimal label $L^*$ representing the column with the maximum reduced cost
\State Initialize $L^* \leftarrow (1, -w_3^*, m^{\text{dry}}, \{1\}, \mathbf{0}_{n_r n_{rv}})$ \label{line:initial}
\State Initialize Priority Queue $\mathcal{Q} \leftarrow \{L^*\}$
\State $\bar{c}^* \leftarrow 0$
\State $\Lambda = \emptyset$
\While{$\mathcal{Q}$ is not empty}
    \State Pop label $L = (i, \phi, u, \boldsymbol{\gamma}, \boldsymbol{\rho})$ with highest $\phi$ from $\mathcal{Q}$
    \State $\phi^U \leftarrow \sum_{j \in S_T \backslash \text{set}(\boldsymbol{\gamma})} \max\{0, p_j - w^*_{1,j}\}$ \Comment{Define upper bound of potential rewards}
    \State $m^p \leftarrow u - m^{\text{dry}} - \sum_{j \in \text{set}(\boldsymbol{\gamma}) \cap S_T} s_j$ \Comment{Calculate required onboard propellant at node $i$}

    \If {${\phi} + \phi^U - \lambda {m^p} \leq \bar{c}^*$} \Comment{If upper bound does not exceed current best cost}
        \State \textbf{continue}
    \EndIf

    \If {$\exists (\phi', u') \in \Lambda_{i, \text{set}(\boldsymbol{\gamma})}$ such that $\phi' \geq \phi$ \textbf{and} $u' \leq u$} \Comment{If the current label is dominated} \label{line:dominate check}
        \State \textbf{continue}
    \EndIf
    \State $\Lambda_{i, \text{set}(\boldsymbol{\gamma})} \leftarrow \{ (\phi', u') \in \Lambda_{i, \text{set}(\boldsymbol{\gamma})} \mid \phi < \phi' \textbf{ or } u > u' \}$ \Comment{Collect non-dominated Pareto front}
    \State $\Lambda_{i, \text{set}({\boldsymbol{\gamma}})} \leftarrow \Lambda_{i, \text{set}({\boldsymbol{\gamma}})} \cup \{ (\phi, u) \}$ \Comment{Add current label to Pareto front}
    
    \If{$i = 0$} \Comment{if $i$ is a starting depot}
        \If{$\omega(\boldsymbol{\gamma}) \in \Omega_c$} \Comment{If the route $\omega$ (corresponds to $\boldsymbol{\gamma}$) already calculated}
        \State \textbf{continue}
        \EndIf
        \State $\bar{c} \leftarrow \phi - \lambda m^p$ \label{line:propellant}
        \If{$\bar{c} > \bar{c}^*$} \Comment{If a higher reduced cost is found}
            \State $\bar{c}^* \leftarrow \bar{c}$
            \State $L^* \leftarrow L$
        \EndIf
        \State \textbf{continue}
    \EndIf
    
    \For{available predecessor node $j$}
        \State $u \leftarrow u / \hat{\mu}_{ji}$ \label{line:u-update}
        \If{$u > m_{\text{max}}$} \Comment{If spacecraft mass exceeds the max mass}
            \State \textbf{continue} 
        \EndIf
        
        \If{$j \in S_T$} \Comment{If the predecessor is a target satellite}
            \If{$s_j + \sum_{i' \in {\text{set}(\boldsymbol{\gamma}) \cap S_T}} s_{i'} > q_{\text{max}}$ \textbf{or} $s_j + u > m_{\text{max}}$} \Comment{If payload or spacecraft max capacity is exceeded}
                \State \textbf{continue} 
            \EndIf
            \State $\phi \leftarrow \phi + (p_j - w_{1,j}^*)$ 
            \label{line:dual profit}
            \State $u \leftarrow u + s_{j}$ \label{line:uj target}
        \ElsIf{$j \in S_R$} \Comment{If the predecessor is a refueling station}
            \State Calculate refuel amount $r_{j}$ from Eq.~\eqref{eq:refueling-rule}
            \State $\phi \leftarrow \phi - (w_{2,{j'(j)}}^* + r_{j}(\lambda + w_{4,{j'(j)}}^*))$ \label{line:refueling-related cost}
            \State $u \leftarrow u - r_{j}$ \label{line:uj refuel}
            \State $\rho_{j} \leftarrow r_j$
        \EndIf
        \State $\boldsymbol{\gamma} \leftarrow [\boldsymbol{\gamma}, j]$
        
        \State Push $L = (j, \phi, u, \boldsymbol{\gamma}, \boldsymbol{\rho})$ to $\mathcal{Q}$
    \EndFor
\EndWhile
\State \textbf{return} $L^*$
\end{algorithmic}
\end{algorithm}

As in the arc-based formulation (Sec.~\ref{Sec:sol-frame-arc}), the trajectory optimization problem for the generated column is also required. Instead of solving $P_l$ in Line~\ref{line:Pl} of Algorithm \ref{Alg:arc-based}, the backward labeling algorithm (Algorithm \ref{Alg:labeling}) is executed to find the visiting sequence and refueling plan with the highest reduced cost. Then, as in Algorithm \ref{Alg:arc-based}, the $P_t$ is solved to sequentially update the $\hat{\mu}$ matrix, and the iteration ends when the $\Delta v$ matrix converges within a predefined threshold or when a cycle is detected.

\subsubsection{Overall Solution Framework of Path-based Formulation}
\begin{figure}[!hbt]
    \centering
    \includegraphics[width=3.25in]{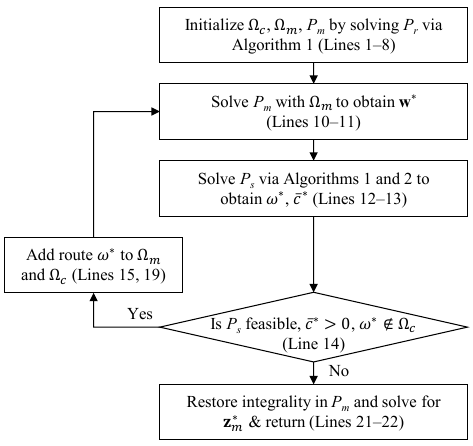}
    \caption{Flow chart of the Algorithm \ref{Alg:cg-framework}}
    \label{fig:flow-chart-alg3}
\end{figure}

The overall solution procedure for the path-based formulation follows a typical column generation framework. Algorithm \ref{Alg:cg-framework} summarizes this procedure and Fig.~\ref{fig:flow-chart-alg3} provides a visual summary. The restricted master Problem ($P_m$) is initialized with routes with a single satellite visit (i.e., solving $P_r$ for each target satellite using Algorithm \ref{Alg:arc-based}, replacing $P_l$ with $P_r$). In each iteration, $P_m$ is solved to obtain the optimal dual variables ($\mathbf{w}^*_1, \mathbf{w}^*_2, w^*_3, \mathbf{w}^*_4$). These dual values are passed to the subproblem ($P_s$), which is solved by executing the Algorithm \ref{Alg:arc-based}, replacing $P_l$ with $P_s$, which uses Algorithm \ref{Alg:labeling} to find the highest reduced cost column. If there is a feasible column that has not been calculated before and whose reduced cost remains positive, the column is added to the set $\Omega_m$, and $P_m$ is resolved. This process repeats until no columns with positive reduced costs remain. After generating all necessary columns, the integrality constraint is restored, and the problem is solved to obtain the optimal solution.
\begin{algorithm}[hbt!]
\caption{Column Generation Framework for Path-based Formulation}
\label{Alg:cg-framework}
\renewcommand{\algorithmicrequire}{\textbf{Input:}}
\renewcommand{\algorithmicensure}{\textbf{Output:}}
\begin{algorithmic}[1]
    \Require Orbital elements of all nodes, problem parameters
    \Ensure Optimal binary selection of routes $\mathbf{z}_m^*$
    \State Initialize restricted column set $\Omega_m \leftarrow \emptyset$, calculated column set $\Omega_c \leftarrow \emptyset$
    \For{each target satellite $i \in S_T$}
        \State Solve $P_r$ with a single-visit route $\omega$ including $i$ using Algorithm \ref{Alg:arc-based} with $P_l \rightarrow P_r$
        \State $\Omega_c \leftarrow \Omega_c \cup \{\omega\}$
        \If{$P_r$ with $\omega$ is feasible}
            \State $\Omega_m \leftarrow \Omega_m \cup \{\omega\}$
        \EndIf
    \EndFor
    \While{$\mathtt{True}$}
        \State Solve $P_m$ with $\Omega_m$
        \State Get optimal dual variables $\mathbf{w}^* = (\mathbf{w}^*_1, \mathbf{w}^*_2, w^*_3, \mathbf{w}^*_4)$
        \State Solve $P_s$ using Algorithm \ref{Alg:arc-based} with $P_l \rightarrow P_s$ and Algorithm \ref{Alg:labeling}
        \State Obtain {optimal route $\omega^*$ and reduced cost $\bar{c}^*$}
        \If{$P_s$ is feasible \textbf{and} ${\bar{c}^*} > 0$ \textbf{and} ${\omega^*} \notin \Omega_c$}
            \State $\Omega_m \leftarrow \Omega_m \cup \{{\omega^*}\}$
        \Else
            \State \textbf{break}
        \EndIf
        \State $\Omega_c \leftarrow \Omega_c \cup \{{\omega^*}\}$
    \EndWhile
    \State Restore integrality constraints in $P_m$ and solve: $\mathbf{z}_m \in \{0, 1\}^{|\Omega_m|}$
    \State \textbf{return} $\mathbf{z}_m^*$
\end{algorithmic}
\end{algorithm}

\section{Case Study} \label{Sec:case_study}

This section presents the case study results of the presented problem formulations and the corresponding solution methodologies. A virtual geosynchronous satellite servicing mission serves as the example for this case study. To simplify the analysis, we assume three key factors in the case study and numerical simulation: 1) all servicing spacecraft are assumed to have identical dry mass and maximum payload capacity, 2) each station has the same maximum refuelable propellant mass, and 3) the time required for servicing and refueling is fixed at two days. As these factors are input parameters of the proposed formulations and methodologies, we can relax the assumptions without further adjustment. The key parameters utilized in the case study and numerical simulations are summarized in Table~\ref{tab:table_param}. The spacecraft specifications are based on the previous studies \cite{shimane_orbital_2024, Choi_2025_OrbitalDepotLocation}. Because this case study focuses on geosynchronous satellites, the canonical distance unit (DU) is defined as the geosynchronous orbit radius, 42,164 km. Consequently, the canonical time unit (TU) is derived using Earth's gravitational parameter ($\mu_e$): $1 \text{ TU} = \sqrt{\text{DU}^3 / \mu_e} = 13,713.4 \text{ sec} \approx 3.809 \text{ hr}$.
Furthermore, since the solution methodology involves solving multiple MILP problems, a computation time limit is established for each MILP (100 s in this study). If the solver cannot find an optimal solution within this time, the best feasible solution is used.

The solution frameworks are implemented in Python programming language with GUROBI optimizer \cite{gurobi2025}, and they run on an Intel Core Ultra 9 285K processor with 64 GB RAM. Parallel computing with 5 workers is used for the numerical experiment, with 4 GUROBI threads solving each MILP.

\begin{table}[hbt!]
    \centering
    \caption{Parameters for the case study and numerical experiments}
    \begin{tabular}{lr}
    \toprule
    Parameter                                                                       & Value/Range \\
    \midrule
    Gravitational parameter of the Earth, $\mu_e$, km\textsuperscript{3}/s\textsuperscript{2}      & $3.986 \times 10^{5}$\\
    Canonical Distance Unit, DU, km                 & 42,164\\
    Canonical Time Unit, TU, hours                    & 3.809\\
    Standard acceleration of gravity, $g_0$, m/s\textsuperscript{2}                 & 9.81\\
    Spacecraft's specific impulse, $I_{sp}$, s                                      & 320\\
    Number of duplicated depots, $n_{dv}$, -                                           & 3\\
    Number of duplicated refueling stations, $n_{rv}$, -                                           & 3\\
    Payload mass to service satellite $i$, $s_i$, kg                                      & [50, 100)\\
    Profit for satellite $i$, $p_i$, - & \{1, 2, 3\}\\
    Spacecraft's dry mass, $m^{\text{dry}}$, kg                        & 500\\
    Maximum mass of spacecraft, $m_{\text{max}}$, kg                                  & 2,000\\
    Maximum mass of payload, $q_{\text{max}}$, kg                                  & 200\\
    Maximum refuel mass at each station, $r_{\text{max}}$, kg                                  & 1,000\\
    Servicing/Refueling time at node $i$, $T^s_i$, days & 2\\
    Weighted factor, $\lambda$, kg\textsuperscript{-1} & $0.0005$\\
    Maximum number of iterations, $L_{\text{max}}$ & 20\\
    Maximum mission makespan, $T_\text{max}$, days & 365\\
    \bottomrule
    \end{tabular}
    \label{tab:table_param} 
\end{table}

\subsection{Sample Servicing Scenario for Geosynchronous Satellites with Refueling Stations}
In this section, we focus on a scenario involving two refueling stations ($n_r = 2$) and five geosynchronous satellites ($n_t = 5$). The orbital elements, the necessary payload mass for each satellite, and the corresponding profit associated with each satellite are detailed in Table \ref{tab:orbital_elements}.

\begin{table}[hbt!]
\centering
\caption{Geocentric Orbital Elements, Payload Requirements, and Profits for Case Study}
\label{tab:orbital_elements}
\begin{tabular}{lccccccccc}
\toprule
Type & Name & $a$ [km] & $e$ & $i$ [deg] & RAAN [deg] & $\omega$ [deg] & $M$ [deg] & \makecell{Payload [kg]} & \makecell{Profit} \\ \midrule
Depot & Depot & 42,164 & 0.00 & 10.00 & 0.00 & 0.00 & 0.00 & - & - \\ \midrule
\multirow{2}{*}{Refuel Station} & Station 1 & 42,164 & 0.00 & 5.00 & 0.00 & 0.00 & 90.00 & - & - \\
 & Station 2 & 42,164 & 0.00 & 5.00 & 0.00 & 0.00 & 270.00 & - & - \\ \midrule
\multirow{5}{*}{Target Satellite} & Target 1 & 42,164 & 0.00 & 3.54 & 298.86 & 0.00 & 316.51 & 60 & 1 \\
 & Target 2 & 42,164 & 0.01 & 1.67 & 170.97 & 0.00 & 100.23 & 90 & 3 \\
 & Target 3 & 42,164 & 0.01 & 3.24 & 248.49 & 0.00 & 262.02 & 50 & 2 \\
 & Target 4 & 42,164 & 0.00 & 9.39 & 178.23 & 0.00 & 258.69 & 70 & 3 \\
 & Target 5 & 42,164 & 0.01 & 6.97 & 342.02 & 0.00 & 234.98 & 80 & 1 \\ \bottomrule
\end{tabular}
\end{table}

Table \ref{tab:case-study_comparison} summarizes the case study results for two proposed methods. The results clearly indicate that both methods yield nearly identical solutions, while the path-based method is more than 100 times faster ({2.779 sec vs 401.093 sec}). Therefore, the solving time of the MILP is significantly greater than that of the sequential NLP or labeling algorithm for column generation. The minor difference in the initial propellant/refueling mass is due to a numerical issue (e.g., convergence tolerance) in the NLP solution. To normalize and assess the performance, we introduce the \textit{gap to ideal profit} ($G_{\text{ideal}}$). The ideal profit ($P_{\text{ideal}}$) represents a hypothetical upper bound where all target satellites are serviced with zero transfer costs, defined as:
\begin{align}
    P_{\text{ideal}} = \sum_{i\in S_T} p_i
\end{align}
The gap to ideal profit is then calculated as:
\begin{align}
    G_{\text{ideal}} = \frac{P_{\text{ideal}} - J}{P_{\text{ideal}}}
\end{align}
where $J$ denotes the objective value obtained from either the arc-based or path-based formulations. It should be noted that $G_{\text{ideal}}$ is a normalized distance to a utopian scenario, not a rigorous mathematical optimality gap. While a smaller $G_{\text{ideal}}$ indicates a solution closer to the ideal scenario, the value may remain high even for the true optimal solution. This is because servicing every target satellite may be physically infeasible under the given mission constraints. Although deriving a rigorous upper bound (e.g., via LP-relaxation) is common in standard routing problems, an LP-relaxation can only yield a local bound for a fixed cost matrix at a specific iteration for each method. Because the cost matrix is also a decision variable in the original problem, computing the true global optimal solution is intractable. In the absence of a tight global bound, $P_{\text{ideal}}$ is introduced as a universally applicable, albeit conservative, upper bound for both arc-based and path-based methods. Consequently, $G_{\text{ideal}}$ serves as a relative comparison metric between two formulations rather than the true optimality gap.

\begin{table}[hbt!]
\centering
\caption{Summary and Comparison of Case Study Results}
\label{tab:case-study_comparison}
\begin{tabular}{lcc}
\toprule
Result Parameter & Arc-based Method & Path-based Method \\ \midrule
Total profit collected & 10.000 & 10.000 \\
Total initial propellant, kg & 876.509 & 876.503 \\
Total refueling mass, kg & 461.431 & 461.367 \\
Number of used vehicles (routes) & 2 & 2 \\
Number of visited targets & 5 & 5 \\
Total computation time, sec & 401.093 & 2.779 \\
Gap to ideal profit, \% & 6.690 & 6.689 \\
\bottomrule
\end{tabular}
\end{table}

The comprehensive sequence of the trajectory, departure and transfer times, $\Delta v$, along with the initial and final masses for each trajectory, is outlined in Tables \ref{tab:arc_results} and \ref{tab:path_results}. Two servicing spacecraft are used to service all satellites to collect profits. For instance, \textit{Servicer 1} departs from the depot and transfers to Target 5 at 0.7 TU with 8.49 TU transfer orbit ($\Delta v = $ 0.66 km/s). After arriving at Target 5, the spacecraft provides in-space servicing for 12.60 TU (2 days). The payload mass (e.g., refueling amount or replacing parts) is delivered to Target 5 (80 kg) and departs for Target 3 to continue the service. Note that after finishing the service at Target 4, the spacecraft visits the Refueling Station 1 for en-route refueling. When the spacecraft arrives at the station, it spends all its propellant and refuels with the amount needed to return to the depot (85.53 kg). This strategy is beneficial because if the spacecraft decides to return to the depot without refueling, it must depart with more propellant, increasing overall propellant consumption. However, if the refueling station is too far away, the propellant required to reach it may be greater. So, the optimizer has to balance between the en-route refueling and departure with more propellant, under the spacecraft's visiting plan and maximum mass constraint.

\begin{table}[hbt!]
\centering
\caption{Detailed Mission Sequence: Arc-based Method}
\label{tab:arc_results}
\resizebox{\textwidth}{!}{
\begin{tabular}{lccccc}
\toprule
Trajectory Segment & Departure [TU] & Transfer [TU] & $\Delta v$ [km/s] & Initial Mass [kg] & Final Mass [kg] \\ \midrule
\textit{Servicer 1} & & & & & \\
Depot $\rightarrow$ Target 5  & 0.70  & 8.49 & 0.66 & 1465.32 & 1187.63 \\
Target 5 $\rightarrow$ Target 3  & 21.79 & 6.91 & 0.59 & 1107.63 & 917.74 \\
Target 3 $\rightarrow$ Target 4  & 41.30 & 7.06 & 0.63 & 867.74  & 709.12 \\
Target 4 $\rightarrow$ Refuel Station 1 & 61.47 & 6.05 & 0.77 & 639.12  & 500.00 \\
Refuel Station 1 $\rightarrow$ Depot & 80.13 & 7.82 & 0.50 & 585.53  & 500.00 \\ \midrule
\textit{Servicer 2} & & & & & \\
Depot $\rightarrow$ Refuel Station 2  & 0.02  & 7.82 & 0.50 & 761.19  & 650.00 \\
Refuel Station 2 $\rightarrow$ Target 2 & 20.44 & 5.02 & 0.38 & 860.77 & 761.93 \\
Target 2 $\rightarrow$ Target 1  & 38.60 & 5.93 & 0.29 & 671.93  & 612.39 \\
Target 1 $\rightarrow$ Refuel Station 2 & 57.13 & 4.84 & 0.31 & 552.39  & 500.00 \\
Refuel Station 2 $\rightarrow$ Depot & 74.58 & 4.10 & 0.90 & 665.13  & 500.00 \\ \bottomrule
\end{tabular}}
\end{table}

\begin{table}[hbt!]
\centering
\caption{Detailed Mission Sequence: Path-based Method}
\label{tab:path_results}
\resizebox{\textwidth}{!}{
\begin{tabular}{lccccc}
\toprule
Trajectory Segment & Departure [TU] & Transfer [TU] & $\Delta v$ [km/s] & Initial Mass [kg] & Final Mass [kg] \\ \midrule
\textit{Servicer 1} & & & & & \\
Depot $\rightarrow$ Target 5         & 0.70  & 8.49 & 0.66 & 1465.31 & 1187.62 \\
Target 5 $\rightarrow$ Target 3      & 21.79 & 6.91 & 0.59 & 1107.62 & 917.73  \\
Target 3 $\rightarrow$ Target 4      & 41.30 & 7.06 & 0.63 & 867.73  & 709.12  \\
Target 4 $\rightarrow$ Refuel Station 1 & 61.47 & 6.05 & 0.77 & 639.12  & 500.00  \\
Refuel Station 1 $\rightarrow$ Depot & 80.13 & 7.82 & 0.50 & 585.53  & 500.00  \\ \midrule
\textit{Servicer 2} & & & & & \\
Depot $\rightarrow$ Refuel Station 2       & 0.02  & 7.81 & 0.50 & 761.19  & 650.00 \\
Refuel Station 2 $\rightarrow$ Target 2   & 20.43 & 5.02 & 0.38 & 860.75  & 761.95 \\
Target 2 $\rightarrow$ Target 1           & 38.60 & 5.93 & 0.29 & 671.95  & 612.41 \\
Target 1 $\rightarrow$ Refuel Station 2   & 57.13 & 4.84 & 0.31 & 552.41  & 500.00 \\
Refuel Station 2 $\rightarrow$ Depot      & 74.57 & 4.10 & 0.90 & 665.09  & 500.00 \\ \bottomrule
\end{tabular}}
\end{table}

\begin{figure}[hbt!]
    \centering
    \subfloat[]{
        \includegraphics[width=0.48\textwidth]{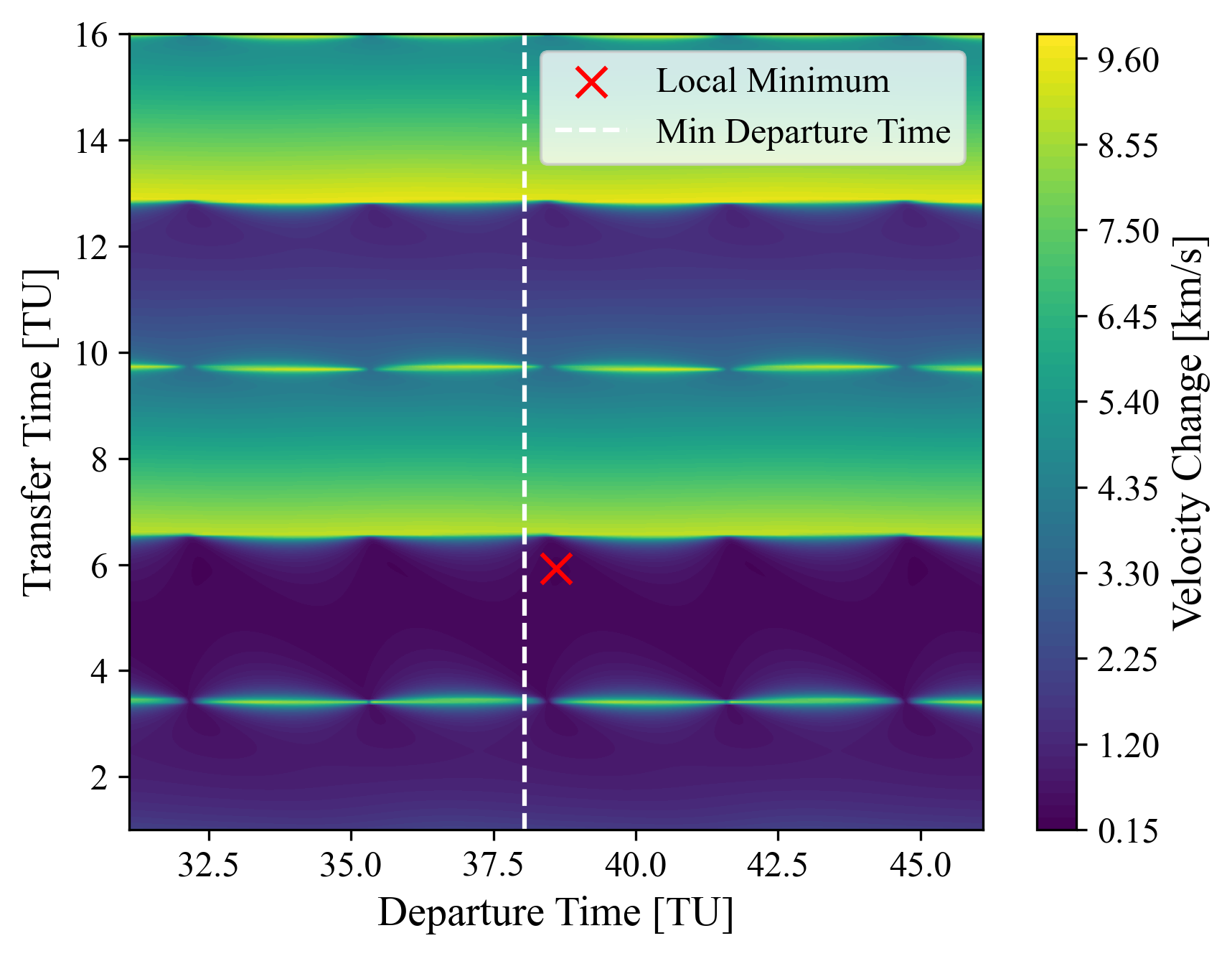} \label{fig:contour_sub1}
    }
    \hfill
    \subfloat[]{
        \includegraphics[width=0.48\textwidth]{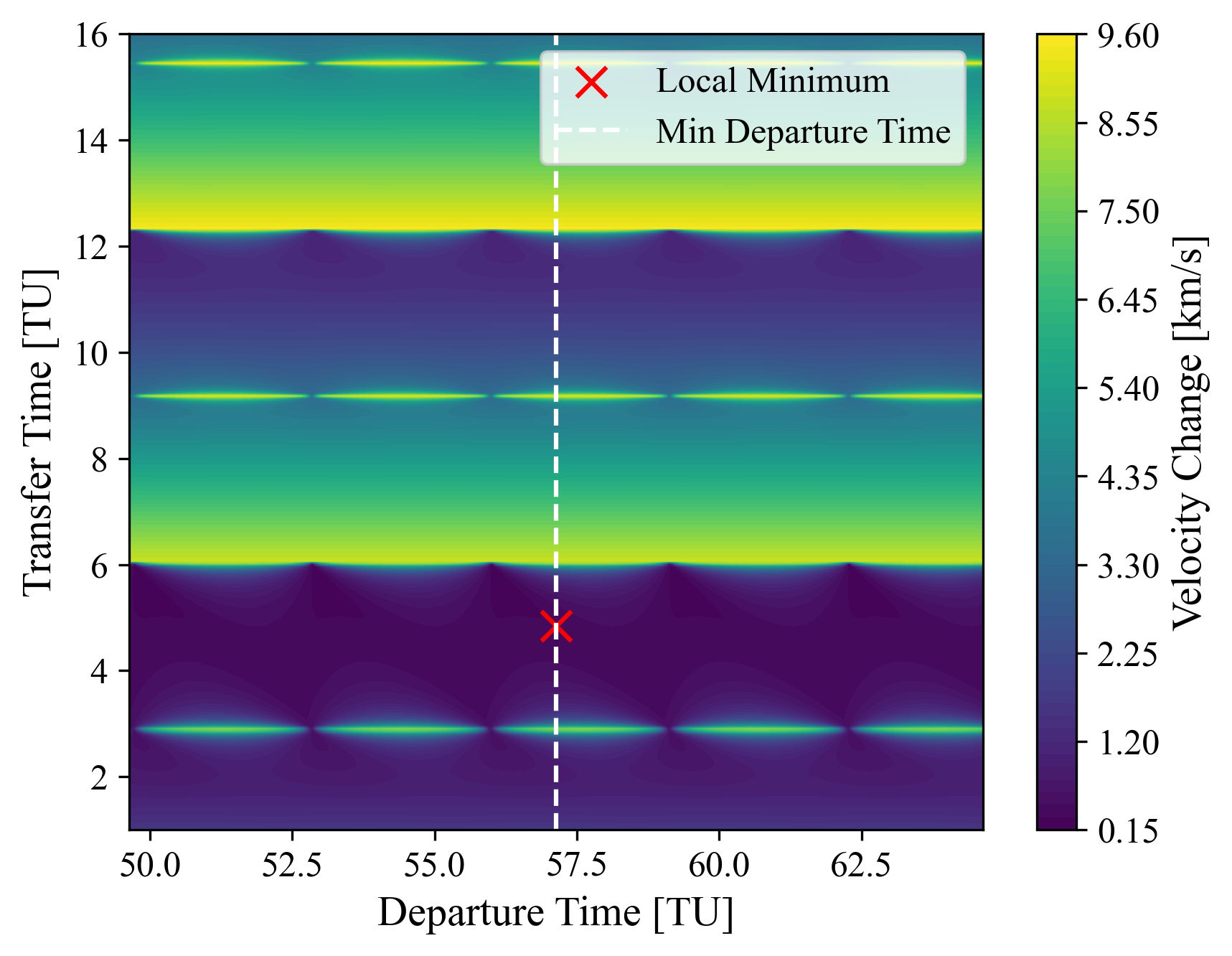} \label{fig:contour_sub2}
    }
    \caption{$\Delta v$ contour plots for two typical transfers: a) Target 2 $\rightarrow$ Target 1, b) Target 1 $\rightarrow$ Refuel Station 2.}
    \label{fig:contour}
\end{figure}

Figure \ref{fig:contour} illustrates the contour plots of $\Delta v$ with different departure and transfer times for two representative trajectory legs in the optimal solution (from Target 2 to Target 1 and from Target 1 to Refuel Station 2 in \textit{Servicer 2} route). The red crosses denote the chosen transfer trajectories, while the white dashed lines represent the minimum departure time for each case. Note that in the case of the trajectory from Target 2 to Target 1, the globally minimum $\Delta v$ trajectory is found (Fig.~\ref{fig:contour_sub1}), while the trajectory from Target 1 to Refuel Station 2 is not a global minimum but a local minimum solution (Fig.~\ref{fig:contour_sub2}) as the SLSQP method is gradient-based. While this trajectory may yield a higher $\Delta v$ compared to the global minimum, it represents a viable trade-off that balances $\Delta v$ with a reduced mission terminal time. If the sole objective is to minimize $\Delta v$, the mission planner might consider alternative nonlinear optimization techniques, such as grid search or metaheuristic approaches, to identify the global optimum.

Additionally, the same servicing scenario with a tight constraint ($T_{\text{max}} = 5$ days) is presented to demonstrate the flexibility of VRTPP-PR. Table \ref{tab:tight-case-comparison} summarizes the results for this scenario. Unlike the previous scenario with a sufficient mission makespan, neither the arc-based nor the path-based method can visit all target satellites. However, both methods yield feasible solutions by abandoning low-profit targets, given the short mission makespan. The path-based method achieves a smaller gap to the ideal profit and requires significantly less computation time (2.636 sec vs. 651.707 sec) in this scenario. While this single instance does not guarantee the universal superiority of the path-based method, it clearly demonstrates that both methods successfully maintain VRTPP-PR's flexible, selective-visiting capability under strict operational limits.

\begin{table}[hbt!]
\centering
\caption{Summary and Comparison of Tight Constraint Scenario Results ($T_{\text{max}} = 5$ days)}
\label{tab:tight-case-comparison}
\begin{tabular}{lcc}
\toprule
Result Parameter & Arc-based Method & Path-based Method \\ \midrule
Total profit collected & 3.000 & 4.000 \\
Total initial propellant, kg & 559.783 & 913.416 \\
Total refueling mass, kg & 0.000 & 0.000 \\
Number of used vehicles (routes) & 1 & 2 \\
Number of visited targets & 1 & 2 \\
Total computation time, sec & 651.707 & 2.636 \\
Gap to ideal profit, \% & 72.799 & 64.567 \\
\bottomrule
\end{tabular}
\end{table}

\subsection{Numerical Experiments}
For scalable comparison of the two solution frameworks, 50 problem instances are generated for each of three refueling station counts ($n_r = 1, 2, 3$) and four target satellite counts ($n_t = 4, 6, 8, 10$). All semimajor axes ($a$) are fixed at 42,164 km (geosynchronous orbit), with eccentricities ($e$) ranging from 0 to 0.01 and inclinations ($i$) varying between 0 and 30 degrees. The right ascension of the ascending node (RAAN), argument of periapsis ($\omega$), and mean anomaly ($M$) can take any value from 0 to 360 degrees. The depot's orbital elements in the case study (see Table \ref{tab:orbital_elements}) are also used in the numerical experiment, except that the inclination was set to 0 deg.

Table \ref{tab:numerical-exp} presents a summary of the result statistics from the numerical experiments. Some instances result in solutions that do not include any target satellites, referred to as trivial problems. These cases are noted in the \textit{\# of Trivial Problems} column. The arc-based formulation is more prone to trivial solutions than the path-based formulation. One of the main reasons is that the arc-based formulation optimizes the overall routing problem simultaneously along with the trajectory optimization problem. If the initial $\Delta v$ matrix is poorly estimated, the optimizer may fail to find a non-trivial solution. Additionally, since each MILP solve is subject to a time limit (100 sec), the solver may fail to find a profitable routing solution to overcome the propellant cost within this limit. On the other hand, the path-based formulation initializes with single-visit columns (routes), enabling the solution procedure to find feasible routes even when the initial trajectory costs are high. 

Additionally, some instances did not converge within 20 iterations, as highlighted in the \textit{\# of Non-converged Problems} column. The computational complexity tends to increase with the number of nodes, leading to more non-converged problems. This scaling trend is also reflected in the mean computation time columns, which become larger as the problem size increases. As in the trivial problem case discussed above, non-convergence issues occur only in the arc-based formulation for this numerical experiment. This is because the path-based method decomposes the problem into smaller subproblems (individual routes) through column generation, making it more robust to convergence failure. Even if the trajectory optimization for a single column (route) does not converge, the other feasible columns remain available in the master problem, ensuring that the overall problem does not fail due to the non-convergence of individual routes.

\begin{table}[htb!]
\centering
\caption{Comparison Results Between Arc-based and Path-based Methods}
\label{tab:numerical-exp}
\begin{tabular}{ccccccccc}
\toprule
\multirow{2}{*}{$n_r$} & \multirow{2}{*}{$n_t$} & \multirow{2}{*}{Formulation} & \multicolumn{2}{c}{Computation Time [sec]} & \multicolumn{2}{c}{\makecell{Gap to\\Ideal Profit [\%]}} & \multirow{2}{*}{\makecell{\# of Trivial\\ Problems}} & \multirow{2}{*}{\makecell{\# of Non-\\Converged\\ Problems}} \\ \cmidrule(lr){4-5} \cmidrule(lr){6-7}
 &  &  & mean & std & mean & std &  &  \\ \midrule
\multirow{8}{*}{1} & \multirow{2}{*}{4}  & Arc-based  & 7.580   & 8.376   & 23.883 & 16.517 & 0  & 0 \\
                   &                     & Path-based & 0.795   & 0.283   & 22.725 & 14.807 & 0  & 0 \\ \cline{2-9}
                   & \multirow{2}{*}{6}  & Arc-based  & 71.528  & 127.783 & 22.757 & 11.695 & 0  & 0 \\
                   &                     & Path-based & 2.206   & 0.676   & 23.244 & 12.056 & 0  & 0 \\ \cline{2-9}
                   & \multirow{2}{*}{8}  & Arc-based  & 336.033 & 324.683 & 24.217 & 8.537  & 0  & 0 \\
                   &                     & Path-based & 3.835   & 1.244   & 25.665 & 9.643  & 0  & 0 \\ \cline{2-9}
                   & \multirow{2}{*}{10} & Arc-based  & 569.843 & 226.748 & 32.259 & 8.107  & 2  & 1 \\
                   &                     & Path-based & 5.542   & 1.613   & 33.832 & 8.608  & 0  & 0 \\ \hline
\multirow{8}{*}{2} & \multirow{2}{*}{4}  & Arc-based  & 445.446 & 332.474 & 18.306 & 12.269 & 0  & 0 \\
                   &                     & Path-based & 2.237   & 0.738   & 18.087 & 12.248 & 0  & 0 \\ \cline{2-9}
                   & \multirow{2}{*}{6}  & Arc-based  & 614.325 & 371.623 & 15.286 & 8.286  & 1  & 0 \\
                   &                     & Path-based & 5.223   & 1.346   & 16.495 & 8.951  & 0  & 0 \\ \cline{2-9}
                   & \multirow{2}{*}{8}  & Arc-based  & 728.035 & 331.824 & 23.320 & 7.665  & 6  & 1 \\
                   &                     & Path-based & 10.489  & 5.077   & 24.365 & 8.276  & 0  & 0 \\ \cline{2-9}
                   & \multirow{2}{*}{10} & Arc-based  & 650.658 & 330.856 & 25.977 & 6.100  & 10 & 0 \\
                   &                     & Path-based & 16.560  & 9.401   & 27.330 & 6.633  & 0  & 0 \\ \hline
\multirow{8}{*}{3} & \multirow{2}{*}{4}  & Arc-based  & 622.581 & 311.380 & 17.467 & 9.100  & 1  & 0 \\
                   &                     & Path-based & 12.821  & 24.815  & 17.724 & 9.311  & 0  & 0 \\ \cline{2-9}
                   & \multirow{2}{*}{6}  & Arc-based  & 731.182 & 369.677 & 14.595 & 6.445  & 2  & 2 \\
                   &                     & Path-based & 25.487  & 24.050  & 15.846 & 6.910  & 0  & 0 \\ \cline{2-9}
                   & \multirow{2}{*}{8}  & Arc-based  & 916.586 & 442.575 & 19.449 & 4.598  & 3  & 2 \\
                   &                     & Path-based & 44.646  & 67.629  & 22.237 & 6.563  & 0  & 0 \\ \cline{2-9}
                   & \multirow{2}{*}{10} & Arc-based  & 760.894 & 384.473 & 27.765 & 6.930  & 3  & 2 \\
                   &                     & Path-based & 88.064  & 84.850  & 29.422 & 8.083  & 0  & 0 \\ \hline
\end{tabular}
\end{table}

\begin{figure}[hbt!]
    \centering
    \subfloat[]{
        \includegraphics[width=3in]{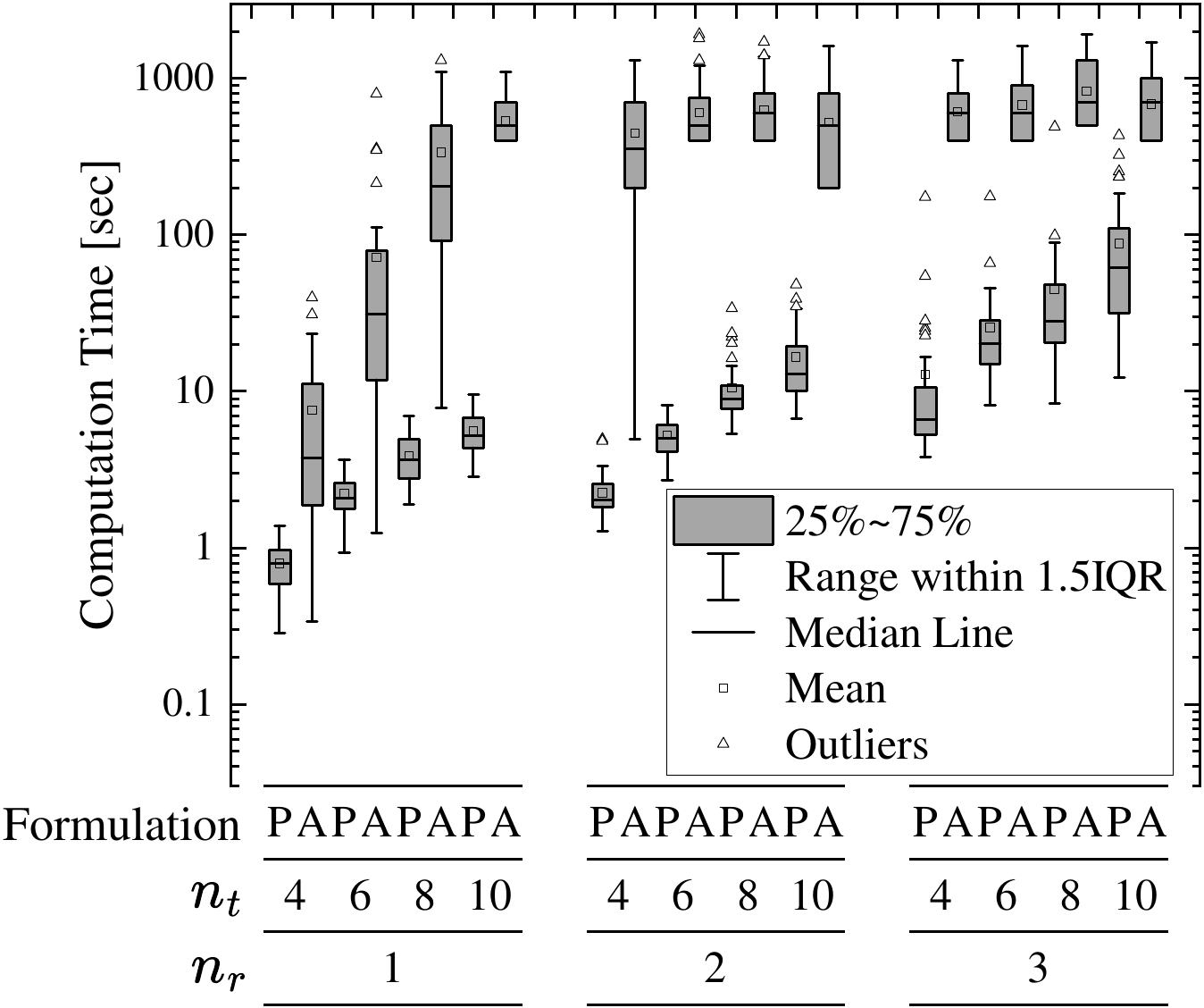}
        \label{fig:num-comp-time}
    }
    \hfill
    \subfloat[]{
        \includegraphics[width=3.2in]{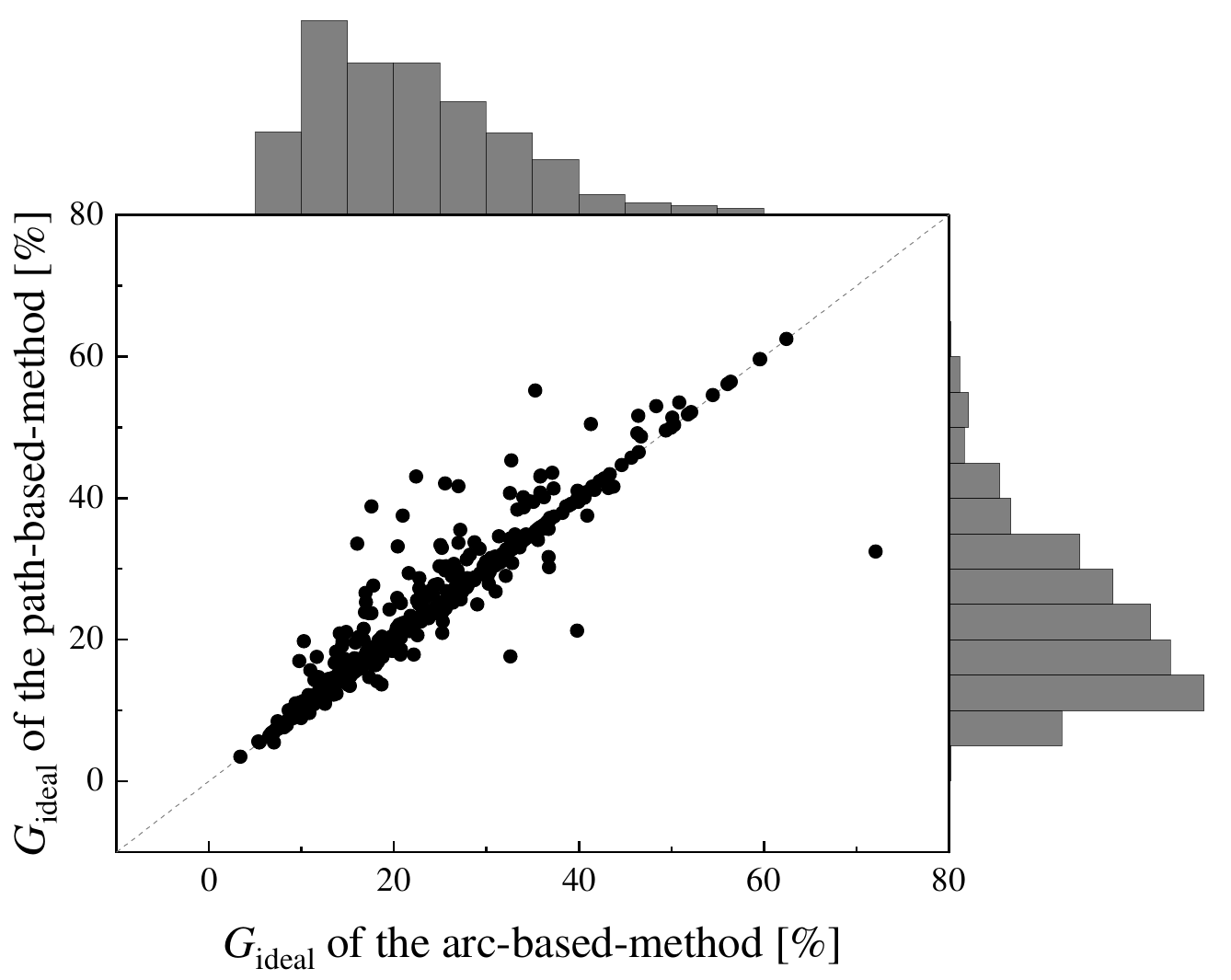}
        \label{fig:g-ideal}
    }
    \caption{Comparison between the arc-based (A) and path-based (P) methods: a) computation time, b) gap to ideal profit.}
    \label{fig:num-exp}
\end{figure}

Figure~\ref{fig:num-exp} visualizes the computation time and the parity plot of the gap to ideal profit for two different formulations, which are also indicated in Table~\ref{tab:numerical-exp}. 
Note that the statistical metrics (i.e., mean and standard deviation) in Table \ref{tab:numerical-exp} and the computation time plots in Figure \ref{fig:num-comp-time} are computed for each method independently, using only instances in which the method converged to a non-trivial solution. Conversely, the parity plots in Figure \ref{fig:g-ideal} are generated exclusively from instances in which both methods converged to non-trivial solutions.
As shown in Fig.~\ref{fig:num-comp-time}, the path-based formulation is significantly faster than the arc-based formulation. Since the maximum computation time of the arc-based formulation is 2,000 sec (20 iterations $\times$ 100 sec limit for each MILP), the indicated computation time gap would be more significant if we do not restrict the maximum computation time. While the path-based formulation's computational time is significantly lower, its objective values are generally lower (worse) than those of the arc-based method, indicating a trade-off between performance and computational time. 

Figure~\ref{fig:g-ideal} shows the comparison result in terms of the gap to ideal profit for both methods. In the figure, the dashed line indicates the $y=x$ line, representing the same $G_{\text{ideal}}$ for both methods. It is worth noting that, depending on the problem, the value of $G_{\text{ideal}}$ can be larger when the problem's constraints are too tight, so that all satellites cannot be visited. However, since the ideal profit value is the same for both formulations, we can still use this parameter as the comparison metric. Since a smaller $G_{\text{ideal}}$ is desired, the cases above the $y=x$ line imply the cases where the arc-based method outperforms the path-based method. Therefore, Fig.~\ref{fig:g-ideal} indicates that most cases show the same or similar results for both formulations, while there is a tendency that the arc-based method outperforms the path-based method.

\subsection{Discussion}
The numerical experiments demonstrate the trade-off between computational time and robustness versus performance (optimality) for the two formulations. Theoretically, the arc-based formulation explores the entire solution space simultaneously. In other words, if the computation time limit and the nonlinear nature of the trajectory costs were ignored, it could eventually identify the global optimum. By contrast, the path-based method cannot provide an optimality guarantee, even under the assumption of linear trajectory costs. Because column generation is designed for the LP, applying it to the LP relaxed MILP and subsequently restoring integrality constraints may cause the algorithm to miss the global optimum. This is because the reduced costs calculated in the subproblem may not accurately reflect the \textit{true reduced cost} of a column once binary constraints are reimposed.

The computational efficiency of the path-based method is highly sensitive to the tightness of spacecraft resource constraints. The performance of the labeling algorithm used in the subproblem depends on the size of the feasible path set ($\Omega_f$). If a spacecraft has sufficient resources (e.g., propellant and payload capacity), the labeling algorithm's search space expands substantially, significantly increasing the time required to generate a single column. On the other hand, the arc-based method exhibits exponential growth in computational complexity as the number of nodes increases, often rendering it intractable on practical time scales. This scaling issue frequently leads the arc-based solver to return only trivial solutions or fail to converge, whereas the path-based method remains more robust; even if a specific subproblem fails to converge, the master problem can still derive a feasible (albeit potentially suboptimal) solution from the existing pool of generated columns.

Depending on the choice of parameters and problem setting, the sensitivity of the objective function to satellite profits can lead to significant variations in the perceived performance gap between the two methods. In this paper, the weighted factor ($\lambda$) is set to be $1/m_{\text{max}}$ to prioritize profit maximization over fuel minimization. Additionally, because the total number of target satellites is relatively small, even missing a single satellite visit can result in a significant gap between the two methods' objective values (e.g., up to 30\% if the ideal profit is 10). Such discrepancies are often driven by the initial conditions of the iterative framework. For instance, the arc-based method may converge to a suboptimal local equilibrium if a \textit{poor} initial $\Delta v$ estimate leads to a suboptimal routing solution in the first iteration. Similarly, the path-based method may exclude potentially optimal columns if the initial cost estimation indicates that the route is infeasible. On the other hand, the arc-based method considers all arcs concurrently; it is more likely to identify higher-profit combinations that the path-based method may have discarded early on.

Note that the trend and comparison of the two formulations discussed in this study are in the context of VRTPP-PR. In this context, the path-based method's superior speed and robustness to non-convergence make it a strong candidate for medium to large-scale constellation servicing with strict resource constraints. However, the potential of the arc-based method to occasionally identify higher-profit solutions by correlating all possible arcs suggests that it may still be preferable for strategic, long-term planning, where computation time is secondary to maximizing scientific or economic profits. There are also space logistics problems which cannot be effectively formulated with the path-based method, such as problems involving not only vehicle routing but also interacting multi-commodity flows. Ultimately, the formulation should be selected based on the mission's specific requirements.

\section{Conclusions} \label{Sec:conclusions}
This paper presents a comparative analysis of arc-based and path-based formulations for the integrated routing and trajectory optimization of satellite servicing missions with partial refueling. By decoupling the interdependent routing, refueling strategy, and transfer trajectory optimizations into a co-optimization framework of MILP and NLP, we addressed the inherent complexity of the VRTPP-PR. Numerical experiments on a geosynchronous servicing scenario demonstrate a fundamental trade-off between the two formulations: the path-based approach using column generation offers superior computational efficiency and robustness, particularly in preventing convergence failures, while the arc-based formulation maintains an advantage in identifying higher-profit solutions by concurrently exploring the entire network.

Future work will extend this framework to multi-depot configurations and missions involving low-thrust electric propulsion, where trajectory sensitivity is more pronounced. Furthermore, integrating spacecraft vehicle design parameters and depot orbital element optimization will provide a more comprehensive tool for long-term, campaign-level mission planning in the evolving space logistics landscape.

\section*{Acknowledgments}
This work was conducted with support from the Air Force Office of Scientific Research (AFOSR), as part of the Space University Research Initiative (SURI), under award number FA9550-23-1-0723. AI technologies were used for grammar checking and refinement. We appreciate Polina Verkhovodova for reviewing and verifying the code.

\bibliography{sample}

\end{document}